\documentclass[a4paper,12pt]{article}
\usepackage[utf8x]{inputenc}
\usepackage{amssymb}
\usepackage{amsmath}
\usepackage{amsthm}
\usepackage{graphicx}
\usepackage{bbm}
\newcommand{\eps}{\varepsilon}

\title{Galerkin-Petrov approach for the Boltzmann equation}
\author{Irene M. Gamba and Sergej Rjasanow}

\begin{document}
\maketitle
\begin{abstract}
In this work, we propose a new Galerkin-Petrov method for the 
numerical solution of the classical spatially homogeneous Boltzmann equation. 
This method is based on an approximation of the distribution function by 
associated Laguerre polynomials and spherical harmonics and test an a 
variational manner with globally defined three-dimensional polynomials. A 
numerical realization of the algorithm is presented.
The algorithmic developments are illustrated 
with the help of several numerical tests.
\end{abstract}
\section{Introduction}
In this paper, we propose a new Galerkin-Petrov method for the 
numerical solution of the classical spatially homogeneous Boltzmann equation. 
This method is based on an approximation of the distribution function by 
associated Laguerre polynomials and spherical harmonics. The test functions 
are polynomials defined globally in $\mathbb{R}^3$. This choice leads to a 
rapid numerical scheme with a high spectral accuracy for smooth solutions. 

Deterministic methods for the Boltzmann equation have been extensively 
studied in the last decades. Overview of these methods can be found, for 
example, in the book of  V. Aristov \cite{ARISTOV:01}  and in a more recent 
review by A. Narayan and A. Kl\"ockner \cite{NARAYAN;KLOECKNER:09}. 
Since the pioneering work of D. Goldstein, B. Sturtevant and J. E. Broadwell 
\cite{GOLDSTEIN;STURTEVANT;BROADWELL:89}, many authors proposed different 
ideas on how to derive a 
discrete version of the Boltzmann collision operator \cite{OHWADA:93},
\cite{ROGIER;SCHNEIDER:94},\cite{WAGNER:95},\cite{PLATKOWSKI;ILLNER:88},
\cite{PALCZEWSKI;SCHNEIDER:98},\cite{PANFEROV;HEINTZ:02}. In 
\cite{KOSUGE;AOKI;TAKATA:01}
the authors studied the difference scheme for a mixture of gases. 
L. Pareschi and G. Russo \cite{PARESCHI;RUSSO:00a},\cite{PARESCHI;RUSSO:00b}
considered deterministic spectral methods for the Boltzmann equation based 
on the Fourier transform. In our paper, we limit our consideration to a 
particular class of deterministic methods, namely, those based on mesh-free 
Galerkin-Petrov discretisation. The main difficulty within the deterministic 
approximation of the Boltzmann collision integral, besides its high 
dimensionality, is the fact that a grid for the integration over the velocity 
space $\mathbb{R}^3$ is not suitable for the integration over the set of all 
directions (i.e. the unit sphere $S^2$). In the case of a regular tensor 
discretisation of the velocity space with $n$ points in each direction, only 
${\cal{O}}(n)$ irregularly distributed integration points would belong to the 
unit sphere. A. Bobylev, A. Palczewski and J. Schneider 
\cite{BOBYLEV;PALCZEWSKI;SCHNEIDER:95} considered this direct approximation 
of the Boltzmann collision integral and showed that the corresponding numerical 
method is consistent. This method requires ${\cal{O}}(n^7)$ arithmetical 
operations per time step and has the formal accuracy of ${\cal{O }}(n^{-1/2})$.
A. Bobylev and S. Rjasanow considered the case of the Maxwell pseudo-molecules 
and utilised an explicit simplification of the Boltzmann equation 
for this model of interaction alongside with the Fast Fourier Transform (FFT) 
to develop a deterministic numerical method \cite{BOBYLEV;RJASANOW:97},
\cite{BOBYLEV;RJASANOW:99}. Their method requires ${\cal{ O}}(n^4)$ arithmetical 
operations per time step and achieves the same low formal accuracy order of 
${\cal{O}}(n^{-1/2})$. A similar method was proposed by L. Pareschi and B. 
Perthame in \cite{PARESCHI;PERTHAME:96}. It appears to be the fastest known 
deterministic numerical method on an uniform grid. At the same time, its 
applications are 
strongly restricted to the case of Maxwell pseudo-molecules. Considering 
the case of hard spheres, A. Bobylev and S. Rjasanow \cite{BOBYLEV;RJASANOW:00}
developed an algorithm, where the integration over the unit sphere is 
completely separated from the integration over the whole space $\mathbb{R}^3$.
The resulting scheme utilises fast evaluation of the generalised Radon and 
X-Ray transforms via the FFT and requires ${\cal{O}}(n^6 \,\log(n))$ operations
per time step with the high formal accuracy of ${\cal{O}}(n^{-2})$. A further 
development of this approach in \cite{GAMBA;HAACK;HAUCK;HU:17}  led to spectral 
schemes for more general collision kernels with a higher efficiency. 
I. Ibragimov
and S. Rjasanow in \cite{IBRAGIMOV;RJASANOW:02}
used a special form of the Boltzmann collision operator, which led to a 
possibility to omit numerical integration over the unit sphere. This idea was 
later used by I. M. Gamba and S. H. Tharkabhushanam 
\cite{GAMBA;THARKABHUSHANAM:09},
\cite{GAMBA;THARKABHUSHANAM:10}, to handle the granular inelastic Boltzmann 
equation. It was developed further in the recent paper \cite{GAMBA;HAACK:14} 
for most general collision cross-section with anisotropic angular scattering 
that includes grazing collisions approximating the Landau collision operator. 
These methods have also been extended to treat systems of Boltzmann equations 
for gas mixtures and multi-energy level gases 
(see \cite{MUNAFO;HAACK;GAMBA;MAGIN:14},
\cite{ZHANG;GAMBA:2-17}). In these extensions of the scheme, the Langrange 
multiplier method is employed to enforce the total conservation properties 
associated with the mixture. The first result on error estimates and 
convergence to Boltzmann-Maxwell equilibrium states for Lagrangian based 
conservative spectral methods for the Boltzmann equation with elastic 
interactions and hard potential with angular cut-off collision kernels was 
published in \cite{ALONSO;GAMBA;THARKABHUSHANAM:17}. A survey of this 
subject can be found in \cite{GAMBA:17}. While the majority of authors use 
an uniform grid in the velocity space, in 
\cite{HEINTZ;KOWALCZYK;GRZHIBOVSKIS:08} A. Heintz, P, Kowalczyk and R. 
Grzhibovskis have used a non-uniform grid.

Reviews of an already substantial amount of publications on the Discrete 
Velocity Models (DVM) for the Boltzmann equation can be found in 
\cite{BERNHOFF;BOBYLEV:07} and in \cite{BOBYLEV;VINERIAN:07}. Constructive 
ideas in this area have been recently proposed by H. Babowsky and his 
co-authors in \cite{BABOVSKY:14a},\cite{BABOVSKY:14b}. Two recent ideas 
regarding the deterministic solution of the Boltzmann equation are the use of 
the Galerkin schemes based on global basis functions (see unpublished 
manuscripts \cite{FONN:13},\cite{KITZLER;SCHROEBERL:13}) and the 
approximation by means of three-dimensional algebraic tensors 
\cite{IBRAGIMOV;RJASANOW:09},\cite{ BEBENDORF;KUENEMUND;RJASANOW:13}. We 
refer to the recent monograph by B. Shizgal \cite{SHIZGAL:15} devoted to the 
spectral methods and an enormous amount of cited literature therein. 
The approach most similar to ours can be found in 
\cite{ENDER;ENDER:99}. Its realization for a rather simple isotropic 
situation is published in \cite{ENDER;ENDER:07}. 

This paper is organised as follows. In Section 2, we give a short description 
of an initial value problem for the Boltzmann equation and present different 
collision kernels. In Section 3, an abstract version of Galerkin-Petrov 
method for a general bilinear operator is formulated. We describe a set of 
basis and test functions in terms of classical polynomials and spherical 
harmonics. Furthermore, the mass and collision matrices are presented in all 
details. A numerical realization of the algorithm is described in Section 4. 
Here, we use a numerical integration for the entries of the mass and 
collision matrices and describe possible time integration schemes. Finally, 
in Section 5, we present the results of numerical computations done by the 
new method for different initial value problems and different collision 
kernels. Conclusions and an outlook can be found in Section 6.
\section{Boltzmann equation}
\label{Sec Introduction}
We consider the initial value problem for the classical 
spatially homogeneous Boltzmann equation
\begin{eqnarray}
\label{Eqn Boltzmann equation}
 \frac{\partial}{\partial t}f(t,v) = Q(f,f)(t,v)\,,\quad 
 t \in \mathbb{R}_+\,,\ v \in \mathbb{R}^3\,,
\end{eqnarray} 
which describes the time evolution of the probability density
\begin{eqnarray}
\label{Eqn Particle density}
 f\ :\ \mathbb{R}_+ \times \mathbb{R}^3 \rightarrow \mathbb{R}_+ \nonumber
\end{eqnarray}
from its initial value 
$$
 f(0,v)=f_0(v)
$$ 
to the final Maxwell distribution
\begin{eqnarray}
\label{Eqn Maxwell distribution}
 \lim_{t \rightarrow \infty} f(t,v) = f_M(v)=\dfrac{\rho_0}{(2\pi\,T_0)^{3/2}}\,
 \mbox{e}^{-\dfrac{|v-V_0|^2}{2\,T_0}}\,.
\end{eqnarray} 
The right-hand side of the equation
\eqref{Eqn Boltzmann equation}, known as the collision
integral or the collision term, has the form
\begin{eqnarray}
\label{Eqn Collision integral}
 Q(f,f)(t,v)=\int\limits_{\mathbb{R}^3}\int\limits_{S^2} B(v,w,e) 
 \Big(f(t,v')f(t,w')-f(t,v)f(t,w) \Big)\,de\,dw\,.
\end{eqnarray} 
Here
$v,w \in \mathbb{R}^3$ are the post-collision velocities,
$e \in S^2 \subset \mathbb{R}^3$ is a unit vector,
$v',w'  \in \mathbb{R}^3$ are the pre-collision velocities, and
$B(v,w,e)$ is the collision kernel.
The operator $Q(f,f)$ represents the change of the distribution function
$f$ due to the binary collisions between particles. A single
collision results in the change of the velocities of the colliding partners
\begin{eqnarray}
\label{Eqn Collision}
 v',w'\ \rightarrow v,w\,. 
\end{eqnarray} 
The  reversible or elastic collision transformation \eqref{Eqn Collision} 
conserves the momentum and the energy
\begin{equation}
\label{Eqn Conservation}
 v+w=v'+w'\,,\quad |v|^2+|w|^2=|v'|^2+|w'|^2\,, 
\end{equation} 
implying that the post- and pre-collisional relative velocities 
$u=v-w$ and $u'=v'-w'$, respectively,  have the same magnitude, i.e. 
$|u'|=|u|$.
The renormalised pre-collisional relative velocity $u'$ defines 
the scattering direction denoted by the unit vector $e$, namely 
$$
e={u'}{|u'|^{-1}}={u'}{|u|^{-1}}\,.
$$
In particular, the conservative exchange 
of binary states  \eqref{Eqn Conservation}  can be written in the following  
centre of mass - relative velocity coordinates form
\begin{eqnarray}
\label{Eqn Collision Transformation}
 v'=\frac{1}{2}\Big(v+w+|u|e\Big)\,,\quad 
 w'=\frac{1}{2}\Big(v+w-|u|e\Big)\,, \quad e \in S^2\,. \nonumber
\end{eqnarray}
In this frame of reference,  the collision kernel, or transition 
probability rate from 
the pre to post states, is, in general, a mapping
\begin{equation}
\label{Eqn coll-ker}
B\ :\ \mathbb{R}^3\times\mathbb{R}^3\times S^2 \rightarrow \mathbb{R}_+.
\end{equation}
It usually is written in a form of a product of a power 
function of the relative speed and a scattering angular function  
\begin{equation}
 \label{Eqn Collision Kernel}
 B(v,w,e)=B\Bigg(|u|,\dfrac{(u,e)}{|u|}\Bigg)=
 C_{\lambda}\,|u|^{\lambda}\,b\Bigg(\dfrac{(u,e)}{|u|}\Bigg)\,,\quad 
 -3 < \lambda \le 1\,.
\end{equation}
These kernels include hard spheres
($\lambda=1$ and $b=1$), hard potentials ($0<\lambda<1$), Maxwell 
pseudo-molecules ($\lambda=0$), and soft potentials 
models ($-3<\lambda<0$). 
In addition, the weak formulation associated to the Boltzmann equation can be 
derived using the binary structure, the conservative collision law, and 
the the symmetries of the collision kernel with respect to the exchange of 
variables (\ref{Eqn Collision Transformation}).
This weak form reads 
\begin{eqnarray}
\label{weak-form}
 && \frac{\partial}{\partial t} \int\limits_{\mathbb{R}^3} f(t,v)\psi(v)\,dv=
 \int\limits_{\mathbb{R}^3} Q(f,f)(t,v) \psi(v)\, dv \nonumber \\ \\
 &=&\!\!\!\! \int\limits_{\mathbb{R}^3}\int\limits_{\mathbb{R}^3} 
 f(t,v)f(t,w) \!\!\int\limits_{S^2} \!\!B(v,w,e)
 \Big(\psi(v')+\psi(w')-\psi(v)-\psi(w)\Big)\,de\,dw\,dv \nonumber
\end{eqnarray} 
for any test function $\psi$ that makes this integral finite. Note that 
in this weak formulation $\psi(v')$ and $\psi(w')$
are the evaluations in the post-collisional velocities. 
This is what subtlety marks the stability of the Boltzmann equation through 
the H-Theorem given below. 
Taking $\psi\in \mathrm{span}\{1, v, |v|^2 \}$ and using the 
elastic exchange of coordinates \eqref{Eqn Conservation}, the following 
conserved quantities are found
  \begin{equation*}
 \frac{\partial}{\partial t} \int\limits_{\mathbb{R}^3} f(t,v)  
 \left(\begin{array}{cr}
         1 \\
         v  \\     
         |v|^2
       \end{array} \right)
  dv\ =\   \int\limits_{\mathbb{R}^3} Q(f,f)(t,v) 
   \left(\begin{array}{cr}
         1 \\
         v  \\     
         |v|^2
       \end{array} \right)
          dv =  \left(\begin{array}{cr}
         0 \\
         0 \\     
       0
       \end{array}\right)\,. 
\end{equation*} 
Thus, the functions from the set  $\{1, v, |v|^2 \}$ are called 
collision invariants.

Finally, we recall the H-theorem that can be obtained by testing with 
$\psi=f(t,\cdot)$. 
If $f\in C^1\Big((0,\infty), L^1(\mathbb{R}^3)\Big)$, then 
\begin{eqnarray*}
 && \frac{\partial}{\partial t} \int\limits_{\mathbb{R}^3} f(t,v)\ln f(t,v)\,dv=
 \int\limits_{\mathbb{R}^3} Q(f,f)(t,v) \ln f(t,v)\,dv=\\
 && -\!\!\!\!\!\!\!\!\!\!\int\limits_{\mathbb{R}^3\times\mathbb{R}^3\times 
S^{2}}\!\!\!\!\!\!\!\!\left(f(\!v'\!,\!t)f(\!w'\!,\!t) - 
f(\!v,\!t)f(\!w,\!t)\right) \dfrac{\ln (f(\!v'\!,\!t)f(\!w'\!,\!t))}{
\ln (f(\!v,\!t)f(\!w,\!t)))} B(v,w,e)de dw dv \leq 0\, .
 \end{eqnarray*}
As anticipated in \eqref{Eqn Maxwell distribution}, the Boltzmann 
H-theorem ensures that the unique stationary equilibrium state is a 
Maxwell distribution, whose moments are the 
the same as those of the initial state.
In addition, this stationary equilibrium state is stable with convergence rates 
depending on the potential rates $\lambda$ and the integrability properties of  
the angular part $b$.
We assume that the angular part  $b$ of the collision kernel is integrable 
over $e\in S^2$.  
If, in addition, the angular function $b$ is bounded, 
this condition is referred as the Grad's cut-off. 
The integrability condition of the angular part $b$ 
implies that the collision operator $Q(f,f)$ splits into a difference of 
two 
positive operators, 
\begin{equation*}
Q(f,f)(t,v)=Q^+(f,f)(t,v) - Q^-(f,f)(t,v)=Q^+(f,f)(t,v)-f(t,v)\,\nu(t,v),
\end{equation*}
where 
\begin{equation*}
Q^+(f,f)(t,v)=\int\limits_{\mathbb{R}^3}\int\limits_{S^2} B(v,w,e) 
 f(t,v')f(t,w')\,de\,dw
\end{equation*}
is the gain operator,
and 
\begin{equation*}
Q^-(f,f)(t,v)=f(t,v)\,\nu(t,v)
\end{equation*}
is the loss operator, provided that the collision frequency 
integral 
\begin{equation*}
\nu(t,v)=\int\limits_{\mathbb{R}^3}\int\limits_{S^2} B(v,w,e) 
 f(t,w)\,de\,dw
\end{equation*}
is well defined.
Without loss of generality, we assume 
\begin{equation}\label{scatt}
 \frac 1{4\pi} \int\limits_{S^2} b\Big(\dfrac{(u,e)}{|u|}\Big)\,de=  
 \frac{1}{2} \int\limits_{-\pi}^{\pi} b(\cos\theta) \sin\theta\, d\theta=1\,.
\end{equation}

It is important to point out, that the case $\lambda=-3$, 
corresponding to the Coulomb interaction, can not be 
modelled by the Boltzmann equation if the function 
$b(\cos\theta)= \cos((u,e){|u|^{-1}}$ is integrable. This is due to the 
divergence of the integral of $f*|u|^{-3}$
in $3$-dimensions for any integrable $f(t,\cdot)$ in 
$v$-space. The loss operator $Q^-(f,f)$ is not well defined in this case.

We will also consider the special forms of isotropic cut-off 
kernel $B$, namely the Variable Hard Spheres model (VHS), see 
\cite{BIRD:81}. In this model the angular dependence of the 
scattering is isotropic, i.e. independent of the scattering angle
\begin{eqnarray}
\label{Eqn Variable Hard Spheres}
 B(v,w,e)=C_{\lambda}\,|u|^{\lambda}\,,\quad -3 < \lambda \le 1\,.
\end{eqnarray}
Our approach can be extended to treat many interesting  non-cut-off 
collision kernels, in which the angular scattering function $b(\cos(\theta))$ 
becomes singular as the scattering angle $\theta$ approaches zero, or 
equivalently 
\begin{equation}
\cos\theta =\dfrac{(u,e)}{|u|}\rightarrow 1\,. 
 \end{equation}
This limit can be associated to a singular behaviour for near grazing 
collisions corresponding to interactions where $v'\approx v$ and $w\approx w$.
Indeed, by the conservative interaction law the relation
\begin{equation}
\label{Eqn Norms Squares}
|v'-v|^2= |w'-w|^2 = |u|^2\,\frac{1-\cos\theta}{2}\,,
\end{equation}
or equivalently  
$$
|v'-v|= |w'-w| = |u| \sin\frac{\theta}{2} 
$$
holds. This implies, that $|v'-v| \approx 0$ is equivalent to  
$\sin(\theta/2)\approx 0$ independently on the norm of the 
relative speed $|u|$.   
 
While we will not cover the non cut-off case in this study. We expect, 
however, that an application of our proposed Galerkin-Petrov scheme will 
address this case as well. It can be done along the lines of the references 
\cite{ZHANG-THESIS:14}, \cite{ZHANG;GAMBA:17}, where a classical 
Discontinuous Galerkin, or a non-conformal Finite Element Method, was 
developed to compute the spectrum of the linearised Boltzmann equation for 
angular non cut-off scattering kernels ranging from hard to soft potentials.
 
The computational approach for the non cut-off case in these studies uses the 
weak formulation \eqref{weak-form}  with the second order Taylor expansion of 
the test function terms $\varphi(v')-\varphi(v)$. This makes it possible to 
perform the cancellation of non-integral angular singularities analytically, 
i.e. by means of the relation (\ref{Eqn Norms Squares}). Thus, a sound numerical
scheme, which is able to handle proper Rayleigh quotients, is formulated. 

A novel way to numerically compute Rayleigh quotients for solutions of the 
linearised 
radial Landau equation by means of Laguerre polynomial expansion  
can be found in a recent publication \cite{BOBYLEV;GAMBA;ZHANG:17}. 
This work relates to our Galerkin Method approach, since it 
indicates, that we can handle the spectral analysis of general, non-radial 
solutions of both the linearised Boltzmann and Landau equations. We will 
elaborate on this feature of the method in an upcoming paper.  

The fast solver derived in this paper can be used to compute anisotropic 
collisions for grazing limits. This allows for obtaining 
approximation rates of the Landau operator by a sequence of Boltzmann operators,
similarly as it was 
done in \cite{GAMBA;HAACK:14}, where a spectral Lagrangian constrains method 
was employed. One starts by solving the initial value problem for the non-linear 
Boltzmann equation \eqref{Eqn Boltzmann equation}-\eqref{Eqn Collision 
integral} in $3$-dimensions in velocity space with the Coulomb interaction 
($\lambda=-3$). 
The collision kernels are given by a $2$-parameter family 
$(\eps,\delta) \in (0,1]\times [0,2)$ of cut-off angular cross sections as
\begin{align}
\label{H-cross-3}
 b_\eps^\delta\Big(\dfrac{(u,e)}{|u|}\Big)=b_\eps^\delta(\cos\theta)=
 -\frac{4}{2\pi H_\delta(\sin(\eps/2))}\,\,
 \frac1{\cos^{3+\delta}\theta}\,\,
 {\mathbbm{1}}_{\cos\theta\geq \sin(\eps/2)}\,.
\end{align}
with
\begin{align}\label{H-cross-4}
H_\delta(x) = 
\begin{cases} 
& \log x \, ,\ \ \text{for}\  \delta=0\,, \\
&-\frac{x^{-\delta}}{\delta},\ \ \text{for} \ 0<\delta<2\,. 
\end{cases}
\end{align}
Note that the case $\delta=0$ corresponds to the Rutherford cross section. 
The corresponding Landau operator limit is independent of  the angular 
scattering cross section $b_\eps^\delta$. Omitting 
the time variable, it can be written as
\begin{equation*}
 Q_L (f,\!f)(v) \!=\! \text{div}_v \Big( \!\int\limits_{\mathbb{R}^3}\!\! 
 |u|^{\lambda+2}\Big(I-\frac{u \otimes u}{|u|^2}\Big) 
 \Big(f(w)\nabla_v f(v)-f(v)\nabla_w f(w)\Big)\,dw\Big)\,.
\end{equation*}
The value $\delta=0$ is the smallest possible exponent when it is possible 
to obtain the Landau equation.
For any value $\delta>2$, however, it is impossible 
to control the higher terms of the expansion (see \cite{GAMBA;HAACK:14}).
This particular case will be the subject of our study an 
upcoming paper.
\section{Galerkin-Petrov approximation}
\label{Sec Galerkin-Petrov approximation}
Let $\mathbb{V}$ be a space of functions with three independent variables and
\begin{eqnarray}
\label{Eqn Bilinear operator}
 Q :\ \mathbb{V}\times\mathbb{V}\rightarrow\mathbb{V}
\end{eqnarray}
a bilinear operator. Let
$$
 f\ :\ \mathbb{R}_+\times\mathbb{R}^3\rightarrow\mathbb{R}
$$
be a time dependent function with
$$
 f(t,\cdot)\in\mathbb{V}\quad \mbox{for all} \; t\in\mathbb{R}_+\,.
$$
We consider an initial value problem
\begin{eqnarray}
\label{Eqn Initial value problem}
 f_t=Q(f,f)\,,\quad \mbox{for} \; t>0\,,\quad f(0,\cdot)=f_0\,. 
\end{eqnarray}
By the use of a finite dimensional subspace $\mathbb{V}_n$ of the 
space $\mathbb{V}$ having a basis
\begin{eqnarray}
\label{Eqn Basis}
 \Phi=\Big(\varphi_1,\dots,\varphi_n\Big)\,,
\end{eqnarray}
we consider an approximation of the function $f$ in the form
\begin{eqnarray}
\label{Eqn Approximation}
 f^{(n)}=\Phi\,\underline{f}=\sum_{j=1}^n f_j\varphi_j\,,\quad 
 \underline{f}\in\mathbb{R}^n\,.
\end{eqnarray}
Furthermore, let 
$$
 \mathbb{V}_n^*\subseteq\mathbb{V}^*
$$
be a finite dimensional subspace of the space $\mathbb{V}^*$ of distributions 
over $\mathbb{V}$ having a basis
\begin{eqnarray}
\label{Eqn Basis^*}
 \Psi=\Big(\psi_1,\dots,\psi_n\Big)\,.
\end{eqnarray}
Then the Galerkin-Petrov scheme for the equation (\ref{Eqn Initial value 
problem}) reads as follows. Find $f^{(n)}(t,\cdot)\in\mathbb{V}_n$ such that
the Galerkin-Petrov equations
\begin{eqnarray}
\label{Eqn Galerkin-Petrov equations}
 \dfrac{d}{dt}<f^{(n)}(t,\cdot),\psi_i>=
 <Q(f^{(n)}(t,\cdot),f^{(n)}(t,\cdot)),\psi_i>\,,\quad
 i=1,\dots,n
\end{eqnarray}
with the initial condition
\begin{eqnarray}
\label{Eqn GP initial condition}
 <f^{(n)}(0,\cdot),\psi_i>=
 <f_0,\psi_i>\,,\quad i=1,\dots,n
\end{eqnarray}
are satisfied for $t>0$. Here, the brackets $<\cdot,\cdot>$ denote the 
action
of the distribution $\psi_i\in\mathbb{V}^*$ on a function from $\mathbb{V}$.
The system (\ref{Eqn Galerkin-Petrov equations}) is in fact a system of 
ordinary differential equations for the time-dependent coefficients $f_j$ of the vector 
$\underline{f}\in\mathbb{R}^n$. By the use of the bilinear structure of the 
operator $Q$, we get a shorter form of the system (\ref{Eqn Galerkin-Petrov 
equations})
\begin{equation}
 \label{Eqn ODE System}
 \dfrac{d}{dt}\Big(M\underline{f}(t)\Big)_i=
 \underline{f}(t)^\top\,Q_i\,\underline{f}(t)\,,\quad
 i=1,\dots,n
\end{equation}
and 
$$
 M\underline{f}(0)=\underline{f}_0\,,\quad 
 \big(\underline{f}_0\big)_i=<f_0,\psi_i>\,,\quad i=1,\dots,n\,.
$$
The matrices $Q_i$ have the entries of the following form
$$
 Q_i[k,\ell]=<Q(\varphi_k,\varphi_\ell),\psi_i>\,,\quad i,k,\ell=1,\dots,n\,,
$$
while the mass matrix $M$ is defined as
$$
 M[i,j]=<\varphi_j,\psi_i>\,,\quad i,j=1,\dots,n\,.
$$
Turning back to the Boltzmann equation, we assume that the initial condition 
$f_0$ belongs to the Schwartz space $\mathbb{S}$ of infinitely smooth functions 
all of whose derivatives are rapidly decreasing. Then the solution $f$ of the 
Boltzmann equation $f(t,\cdot)$ is again a Schwartz space function for all 
times $t$, see \cite{DUDUCHAVA;RJASANOW:05}. Thus, the basis functions 
$\varphi_j$ belong to the subspace
$$
 \mathbb{S}_n=\mbox{span}\Phi\subset \mathbb{S}\,.
$$
The dual space $\mathbb{S}^*$ is the space of tempered 
distributions. The space $\mathbb{S}^*$ contains among others polynomials 
of arbitrary degree.
\subsection{Basis functions}
In this subsection, we introduce a set of globally defined basis functions.
\subsubsection{Classical polynomials and spherical harmonics}
First, we give the definitions and the main properties of the associated 
Laguerre polynomials, associated Legendre polynomials, and of the spherical 
harmonics.
\subsubsection*{Associated Laguerre polynomials}
The classical associated Laguerre polynomial of degree $k$
is the polynomial solution of the differential equation
$$
 x\,y''+(\alpha-1+x)\,y'+k\,y=0\,,\quad \alpha\in\mathbb{R}_+.
$$
It is denoted by $L_k^{(\alpha)}$. By the use of the abbreviation
$$
 \left(\begin{array}{c} 
        k+\alpha \\ m
       \end{array}
 \right)=\dfrac{(k+\alpha)(k-1+\alpha)\dots(k-m+\alpha)}{m!}\,,
$$
an explicit formula for the polynomial $L_k^{(\alpha)}$ reads
$$
 L_k^{(\alpha)}(x)=\sum_{i=0}^k(-1)^i\left(\begin{array}{c} 
 k+\alpha \\ k-i \end{array}\right)\dfrac{x^i}{i!}\,.
$$
The orthogonality property of the associated Laguerre polynomials can be 
written as
$$
 \int\limits_0^\infty x^\alpha \mbox{e}^{ -x} 
 L_k^{(\alpha)}(x)L_m^{(\alpha)}(x)\,dx=
 \dfrac{\Gamma(k+1+\alpha)}{k!}\,\delta_{k,m}\,,
$$
where $\delta_{k,m}$ is the Kronecker symbol. Thus, the polynomials are 
orthogonal with respect to the measure $x^\alpha \mbox{e}^{-x}\,dx$.
For numerical computations of the associated Laguerre polynomials, we use the 
initial functions
$$
 L_0^{(\alpha)}(x)=1\,,\quad L_1^{(\alpha)}(x)=1+\alpha-x
$$
and the following recursion for $k\geq 2$
$$
 L_k^{(\alpha)}(x) = \frac{(2k-1+\alpha-x)L_{k-1}^{(\alpha)}(x)-
 (k-1+\alpha)L_{k-2}^{(\alpha)}(x)}{k}\,. 
$$
\subsubsection*{Associated Legendre polynomials}
The classical associated  Legendre polynomial
is the polynomial solution of the differential equation
$$
 (1-x^2)\,y''-2x\,y'+\Big(\ell(\ell+1)-\frac{m^2}{1-x^2}\Big)\,y=0\,,
$$
where the index $\ell$ is the degree and $m$ the order of the associated 
Legendre polynomial $P_{\ell,m}$.
An explicit formula for the polynomial $P_{\ell,m}$ is
$$
 P_{\ell,m}(x)=\frac{(-1)^m}{2^\ell\,\ell!}\,(1-x^2)^{m/2}
 \frac{d^{\ell+m}}{dx^{\ell+m}}(x^2-1)^\ell\,,\quad 0\le m\le\ell\,.
$$
The orthogonality properties of the associated Legendre polynomials read as
$$
 \int\limits_{-1}^1P_{\ell_1,m}(x)P_{\ell_2,m}(x)\,dx=
 2\frac{(\ell+m)!}{(2\ell+1)(\ell-m)!}\,\delta_{\ell_1,\ell_2}\,,
$$
for fixed $m$ and, in the case $\ell_1=\ell_2=\ell$. Furthermore,
$$
 \int\limits_{-1}^1\frac{1}{1-x^2}P_{\ell,m}(x)P_{\ell,k}(x)\,dx=
 \left\{\begin{array}{ccc}
         0 & \mbox{for} &m\ne k \\
         \dfrac{(\ell+m)!}{m(\ell-m)!} & \mbox{for} & k=m\neq 0
       \end{array}\right.
$$
for a fixed $\ell$. For $k=m=0$, the last integral diverges.
For numerical evaluations of the associated Legendre polynomials, we use 
the initial functions
$$
 P_{m,m}(x)=(-1)^m (2m-1)!!\, (1- x^2)^{m/2}\,,\
 P_{m+1,m}(x)=x\,(2m+1) P_{m,m}(x)
$$
and the following recursion for $k=m+2,\dots,\ell$
$$
 P_{k,m}(x)=\frac{(2k-1)xP_{k-1,m}(x)-(k-1+m)P_{k-2,m}(x)}{k-m}\,.
$$
\subsubsection*{Spherical harmonics}
The spherical harmonics $Y_{\ell,m}$ are the complete and orthonormal set of 
eigenfunctions of the angular part of the three-dimensional Laplace's equation
$$
 \left(\dfrac{\partial^2}{\partial \theta^2}+
       \dfrac{\cos\theta}{\sin\theta}\dfrac{\partial}{\partial\theta}+
       \dfrac{1}{\sin^2\theta}\dfrac{\partial^2}{\partial \phi^2}\right)
       Y_{\ell,m}(\phi,\theta)=-\ell(\ell+1)Y_{\ell,m}(\phi,\theta)\,,
$$
for $\ell\in\mathbb{N}_0$ and $m=-\ell,\dots,0,\dots,\ell$.
An explicit formula for the spherical harmonics with the 
parameterisation
\begin{equation}
 \label{Eqn Sphere Parameterisation}
 e=\left(\begin{array}{c}
         \cos\phi\,\sin\theta \\ \sin\phi\,\sin\theta \\ \cos\theta
         \end{array}
   \right)
\end{equation}
is
$$
 Y_{\ell,m}(\phi,\theta)=
 \sqrt{\dfrac{2\ell+1}{4\,\pi}\dfrac{(\ell-m)!}{(\ell+m)!}\,}\,P_{\ell,m}
 (\cos\theta)\,\mbox{e}^{ \imath\,m\phi}\,.
$$
Here, $P_{\ell,m}$ are the associated Legendre polynomials.
The orthogonality property of the spherical harmonics reads as
$$
 \int\limits_{S^2}Y_{\ell_1,m_1}(e)Y_{\ell_2,m_2}(e)\,de=
 \delta_{\ell_1,\ell_2}\delta_{m_1,m_2}\,.
$$
However, for our purposes, we will use the real valued version of the spherical 
harmonics in the form
$$
 Y_{\ell,m}(\phi,\theta)=
 \sqrt{\dfrac{2\ell+1}{2\,\pi}\dfrac{(\ell-m)!}{(\ell+m)!}\,}\,P_{\ell,m}
 (\cos\theta)\,\cos(m\phi)
$$
for $m>0$,
$$
 Y_{\ell,0}(\phi,\theta)=
 \sqrt{\dfrac{2\ell+1}{4\,\pi}\,}\,P_{\ell,0}(\cos\theta)
$$
for $m=0$ and
$$
 Y_{\ell,m}(\phi,\theta)=
 \sqrt{\dfrac{2\ell+1}{2\,\pi}\dfrac{(\ell-m)!}{(\ell+m)!}\,}\,P_{\ell,m}
 (\cos\theta)\,\sin(-m\phi)
$$
for $m<0$.
\subsubsection{Basis functions}
In three dimensional spherical coordinates
$$
 v=\varrho\,e_v\,,\quad 0\le\rho<\infty\,,\ e_v\in S^2\,,
$$
we decompose the basis function $\varphi_j$ as follows
$$ 
\varphi_j(v)=\varphi_j(\varrho\,e_v)=\Phi_{k,\ell}(\varrho)\,Y_{\ell,m}(e_v)\,,
\quad
 k\in\mathbb{N}_0\,,\ \ell\in\mathbb{N}_0\,,\ -\ell\le m\le\ell\,.
$$
Thus, the global index $j$ is a function of three indices $j=(k,\ell,m)$. 
Since the angular part of the function $\varphi_j$ is already defined, we 
look at 
the radial part and write the function $\Phi_k$ in the form
$$
 \Phi_{k,\ell}(\varrho)=
 \mu_{k,\ell}\,\mbox{e}^{-\varrho^2/2}\,
 L_k^{(\ell+1/2)}(\varrho^2)\,\varrho^\ell.
$$
The normalisation parameters $\mu_{k,\ell}$ are chosen so, that 
the functions $\Phi_{k,\ell}$ will compose an orthonormal system with respect 
to the 
measure $\varrho^2\,d\varrho$. Setting 
$\varrho^2=x\,,\ 2\,\varrho\,d\varrho=dx$, we get 
\begin{eqnarray*}
\label{Eqn Orthonormal system}
 && \int\limits_0^\infty\mu_{k_1,\ell}\mu_{k_2,\ell}\,\varrho^{2\ell+2}\, 
\mbox{e}^{-\varrho^2}L_{k_1}^{(\ell+1/2)}(\varrho^2)L_{k_2}^{(\ell+1/2)}
(\varrho^2)\,d\varrho= \\
 && 
\frac{1}{2}\,\mu_{k_1,\ell}\mu_{k_2,\ell}\int\limits_0^\infty\varrho^{2\ell+1}
\mbox{e}^{-\varrho^2}L_{k_1}^{(\ell+1/2)}(\varrho^2)L_{k_2}^{(\ell+1/2)}
(\varrho^2)\,2\varrho\,d\varrho= \\
 && \frac{1}{2}\mu_{k_1,\ell}\mu_{k_2,\ell}\int\limits_0^\infty
 x^{\ell+1/2}\mbox{e}^{-x}L_{k_1}^{(\ell+1/2)}(x)L_{k_2}^{(\ell+1/2)}(x)\,dx= \\
 && \frac{1}{2}\,\mu_{k,\ell}^2
 \dfrac{\Gamma(k+\ell+3/2)}{k!}\,\delta_{k_1,k_2}\,,
\end{eqnarray*}
in the case $k_1=k_2=k$. To obtain an orthonormal system, we set
$$
 \mu_{k,\ell}=\sqrt{\frac{2\,k!}{\Gamma(k+\ell+3/2)}\,}\,.
$$
This yields the form of the function $f^{(n)}$ in spherical coordinates 
$v=\varrho\,e_v$ 
$$
 f^{(n)}(v)=\sum_{k=0}^K\sum_{\ell=0}^L\sum_{m=-\ell}^\ell
 f_{k,\ell,m}\Phi_{k,\ell}(\varrho)\,Y_{\ell,m}(e_v)\,.
$$
The number of the basis functions is 
$$
 n=(K+1)\,(L+1)^2\,.
$$
\subsection{Test functions}
All basis functions belong to the Schwartz space $\mathbb{S}$ of infinitely 
smooth functions all of whose derivatives are rapidly decreasing. Thus the 
collision integral $Q(\varphi_k,\varphi_\ell)$ is a Schwartz function as 
well and, therefore, any tempered distribution can be chosen as a test 
function $\psi_i$.  In the case of regular distribution $\psi_i$ 
identified with a continuous function
$$
 \psi_i\,:\,\mathbb{R}^3\rightarrow\mathbb{R}\,,
$$
the entries of the matrices 
$Q_i$ can be 
evaluated for $i,k,\ell=1,\dots,n$ as follows
\begin{eqnarray}
 &&Q_i[k,\ell]=<Q(\varphi_k,\varphi_\ell),\psi_i>
 =\int\limits_{\mathbb{R}^3}Q(\varphi_k,\varphi_\ell)(v)\psi_i(v)\,dv 
 \nonumber \\ 
 &=&\frac{1}{2}\int\limits_{\mathbb{R}^3}\!\int\limits_{\mathbb{R}^3}\!
 \varphi_k(v)\varphi_\ell(w)\int\limits_{S^2}\!B(v,w,e)
\Big(\!\psi_i(v')\!+\!\psi_i(w')\!-\!\psi_i(v)\!-\!\psi_i(w)\!\Big)de\,dw\,dv,
 \nonumber
\end{eqnarray}
where the weak form of the collision integral (\ref{weak-form}) has been 
used. 
If the set of test functions contains a collision invariant, the 
corresponding 
matrices $Q_i$ will vanish completely, and, the corresponding macroscopic 
quantity will be conserved automatically.

One possible choice is a pure Galerkin method with
$$
 \psi_i=\varphi_i\,,\quad i=1,\dots,n\,.
$$
In this case, the mass matrix $M$ is the identity matrix due to the 
orthogonality of the system. However, an additional numerical conservation 
procedure is necessary.

Due to an automatic fulfilment of the conservation properties, the 
following choice of test functions for a index $i=(k,\ell,m)$ seems to be 
natural
$$
 \psi_i(v)=L_k^{(\ell+1/2)}(\varrho^2)\,\varrho^\ell\,Y_{\ell,m}(e_v)\,,\quad 
\mbox{for} \; v=\varrho\,e_v\,.
$$
These globally defined polynomials are in fact the basis functions 
without the factor $\mu_{k,\ell}\,\mbox{e}^{-\varrho^2/2}$.
All 
five collision invariants are included in the set of the test functions, 
namely
\begin{eqnarray*}
 \psi_{0,0,0}(v)&=&\sqrt{\frac{1}{4\,\pi}}\,,\\
 \psi_{0,1,-1}(v)&=&\sqrt{\frac{3}{4\,\pi}}\,\varrho\,\sin\phi\,\sin\theta=
 \sqrt{\frac{3}{4\,\pi}}\,v_2\,, \\
 \psi_{0,1,0}(v)&=&\sqrt{\frac{3}{4\,\pi}}\,\varrho\,\cos\theta=
 \sqrt{\frac{3}{4\,\pi}}\,v_3\,, \\
 \psi_{0,1,1}(v)&=&\sqrt{\frac{3}{4\,\pi}}\,\varrho\,\cos\phi\,\sin\theta=
 \sqrt{\frac{3}{4\,\pi}}\,v_1\,, \\
 \psi_{1,0,0}(v)&=&\sqrt{\frac{1}{4\,\pi}}\,\Big(-\varrho^2+\frac{3}{2}\Big)=
 \sqrt{\frac{1}{4\,\pi}}\,\Big(-|v|^2+\frac{3}{2}\Big)\,.
\end{eqnarray*}
Thus, the conservation properties are now ensured automatically.

For a regular Galerkin-Petrov scheme, it is necessary to choose the same number 
of basis and test functions, i.e. for
$$
 k=0,\dots,K\,,\ \ell=0,\dots,L\,,\ -\ell\le m\le\ell\,,
$$
we get
$$
 n=(K+1)\,(L+1)^2\,.
$$
\subsection{Mass matrix}
The mass matrix $M\in\mathbb{R}^{n\times n}$ has the entries
$$
 M[i,j]=<\varphi_j,\psi_i>\,,\quad i,j=1,\dots,n\,,
$$
where $i=(k_i,\ell_i,m_i)$ and $j=(k_j,\ell_j,m_j)$. Since both, basis and test 
functions contain spherical harmonics which are mutually orthonormal, in 
spherical coordinates we obtain 
$$
M[i,j]=0\,,\quad \mbox{for} \; \ell_i\neq\ell_j\ \mbox{or} \; m_i\neq m_j
$$
and
\begin{equation}
\label{Eqn Mass Matrix Entries} 
 M[i,j]=\mu_{k_j,\ell}\int\limits_0^\infty\varrho^{2\ell+2}
 \mbox{e}^{-\varrho^2/2}\,L_{k_j}^{(\ell+1/2)}(\varrho^2)\,
 L_{k_i}^{(\ell+1/2)}(\varrho^2)\,d\varrho 
\end{equation}
for $\ell_i=\ell_j=\ell$ and $m_i=m_j=m$. Thus, the mass matrix is rather 
sparse and, since $M[i,j]$ do not depend on $m$, has many equal 
non-zero entries.
\subsection{Collision matrices}
For a general interaction model, the collision matrices $Q_i$ have the 
entries
\begin{eqnarray}
 \label{Eqn Collision Matrices Entries}
 Q_i[k,\ell]=\int\limits_{\mathbb{R}^3}\int\limits_{\mathbb{R}^3}
 \varphi_k(v)\varphi_\ell(w)q_i(v,w)\,dw\,dv\,,
\end{eqnarray}
where
\begin{equation}
 \label{Eqn Collision Distribution}
 q_i(v,w)=\int\limits_{S^2}B(v,w,e)\,
 \Big(\psi_i(v')+\psi_i(w')-\psi_i(v)-\psi_i(w)\Big)\,de\,.
\end{equation}
The integration (\ref{Eqn Collision Distribution}) is
an important part of the generation of the collision matrices. For the VHS 
model of interaction (\ref{Eqn Variable Hard Spheres}), and for general 
polynomial test functions, this integration can be done analytically 
leading to a function $q_i$ which is a polynomial in six variables 
$v$ and $w$ multiplied by $|u|^\lambda$.
For more general models of interaction and for test functions given in 
spherical coordinates an analytic integration seems to be impossible.
Furthermore, for  the non cut-off collision models, the kernel $B$ has a 
singularity
and the corresponding numerical integration should be done very carefully.
\section{Numerical realization}
The main advantage of the above Galerkin-Petrov method is the possibility to 
precompute and to store all the collision matrices $Q_i$ and the mass matrix 
$M$ for different discretisation parameters $K$ and $L$. Furthermore, these 
matrices are also independent of  a time discretisation scheme and 
corresponding 
time discretisation parameters. Therefore, once computed, experiments with 
different time discretisation schemes  can be easily  performed.
For the numerical integration over $\mathbb{R}^3$, we will use spherical 
coordinates and a combination of the radial Gauss-Laguerre quadratures with the 
Lebedev quadratures for the integration over the unit sphere.  
For a given function $g\,:\,\mathbb{R_+}\rightarrow\mathbb{R}$, the 
Gauss-Laguerre quadrature is applied to the integrals of the form
$$
 I[g]=\int\limits_0^\infty x^{1/2} \mbox{e}^{-x}\,g(x)\,dx
$$
and results in an approximation
$$
 I_{N_{GL}}[g]=\sum_{i=1}^{N_{GL}}\omega_i^{GL}\,g(x_i)\,.
$$
The weights $\omega_i^{GL}$ and the positions $x_i$ are available for any 
$N_{GL}$ with an arbitrary accuracy (see \cite{SALZER;ZUCKER:49}).
By the use of the parameterisation (\ref{Eqn Sphere Parameterisation}), the 
integral over the unit sphere
for a given function $g\,:\,S^2\rightarrow\mathbb{R}$
$$
 I[g]=\int\limits_{S^2}\,g(e)\,de
$$
can be transformed into the corresponding integrals over the rectangular 
domain 
$[0,2\,\pi]\times[0,\pi]$ and subject to the subsequent application of 
the classical Gauss 
quadratures. However, the Lebedev quadratures 
\cite{LEBEDEV:76},\cite{LEBEDEV:77}
$$
 I[g]=4\,\pi\sum_{j=1}^{N_L}\omega_j^L\,g(e_j)\,,\quad e_j\in S^2
$$
are invariant under finite rotation groups and available for many
values of $N_L$. The first of them are for $N_L=6,14,26,38,50,74,86,110$.
We claim that this set will be sufficient for our first tests.
\subsubsection*{Mass matrix}
The Gauss-Laguerre quadratures will be used for numerical computation of 
the mass matrix entries corresponding to (\ref{Eqn Mass Matrix Entries}). 
With the 
substitution $\varrho^2=x\,,\ 2\,\varrho\,d\varrho=dx$, we get 
$$ 
 M[i,j]=\frac{1}{2}\,\mu_{k_j,\ell}\int\limits_0^\infty x^{\ell+1/2} 
 \mbox{e}^{-x}\Big(x^\ell\mbox{e}^{x/2}\,L_{k_j}^{(\ell+1/2)}(x)\,
 L_{k_i}^{(\ell+1/2)}(x)\Big)\,dx 
$$
and approximate these entries as
$$ 
 M_{N_{GL}}[i,j]=\frac{1}{2}\,\mu_{k_j,\ell}\
 \sum_{i_v=1}^{N_{GL}}\omega_{i_v}^{GL}\,
 x_{i_v}^\ell\mbox{e}^{x_{i_v}/2}\,L_{k_j}^{(\ell+1/2)}(x_{i_v})\,
 L_{k_i}^{(\ell+1/2)}(x_{i_v})\,.
$$
As we have mentioned before, only few entries of the mass matrix are different 
from zero and they are computed numerically during the initialisation. This 
requires just a few seconds of computer time. Then we use LAPACK package to 
perform the LU decomposition of the matrix $M_{N_{GL}}$ in order to solve the 
systems of linear equations with the mass matrix in initial and later in every 
time step of the algorithm. Formally, the numerical work for this decomposition 
is ${\cal{O}}(n^3)$. However, the corresponding computer time is negligible 
in our experiments.
\subsubsection*{Collision matrices}
The computation of the collision matrices is the most important and numerically 
difficult step of the algorithm. However, it is an initialisation step and will 
be done only once for the given collision kernel and for the fixed parameters 
$K,L,N_{GL}$ and $N_L$. Then all $n$ 
collision matrices of the dimension $n\times n$ will be stored and used for 
all computations on the later stages. 
By the use of the 
substitution $\varrho^2=x\,,\ 2\,\varrho\,d\varrho=dx$ again, we get for a 
function $g\,:\,\mathbb{R}^3\rightarrow \mathbb{R}$
\begin{eqnarray*}
 I[g]=\int\limits_{\mathbb{R}^3}g(v)\,dv&=&\int\limits_0^\infty\varrho^2
 \int\limits_{S^2}g(\varrho\,e)\,de\,d\varrho \\
 &=& \frac{1}{2}\int\limits_0^\infty x^{1/2} \mbox{e}^{-x}
 \Big(\mbox{e}^x\int\limits_{S^2}g(\sqrt{x}\,e)\,de\Big)\,dx
\end{eqnarray*}
and approximate these integrals as follows
$$
 I_{N_{GL},N_L}[g]=2\,\pi\,\sum_{i_v=1}^{N_{GL}}\omega_{i_v}^{GL}
 \mbox{e}^{x_{i_v}}
 \sum_{j_v=1}^{N_L}\omega_{j_v}^L g(\sqrt{x_{i_v}}\,e_{j_v})\,.
$$
Thus, for the entries of the collision matrices $Q_l[k,\ell]$ with 
$k=(k_v,\ell_v,m_v)$ and $\ell=(k_w,\ell_w,m_w)$ we get 
\begin{eqnarray*}
  \Big(Q_{N_{GL},N_L}\Big)_i[k,\ell]=
 (2\,\pi)^2&&\hspace{-0.5cm}
 \sum_{i_v=1}^{N_{GL}}\omega_{i_v}^{GL}
 x_{i_v}^{\ell_v}\mbox{e}^{x_{i_v}/2}\,L_{k_v}^{(\ell_v+1/2)}(x_{i_v})\times \\
 &&\hspace{-0.5cm}
 \sum_{j_v=1}^{N_L}\omega_{j_v}^LY_{\ell_v,m_v}(e_{j_v})\times \\
 &&\hspace{-0.5cm}
 \sum_{i_w=1}^{N_{GL}}\omega_{i_w}^{GL}
 x_{i_w}^{\ell_v}\mbox{e}^{x_{i_w}/2}\,L_{k_w}^{(\ell_w+1/2)}(x_{i_w})\times \\
 &&\hspace{-0.5cm}
 \sum_{j_w=1}^{N_L}\omega_{j_w}^LY_{\ell_w,m_w}(e_{j_w})\,
 \big(q_L\big)_i(v_{i_v,j_v},w_{i_w,j_w})\,,
\end{eqnarray*}
where $v_{i_v,j_v}=\sqrt{x_{i_v}}e_{j_v},\,w_{i_w,j_w}=\sqrt{x_{i_w}}e_{j_w}$ 
and
\begin{eqnarray*}
 &&\big(q_L\big)_i(v_{i_v,j_v},w_{i_w,j_w})=
 \sum_{j=1}^{N_L}\omega_{j}^L
 B(v_{i_v,j_v},w_{i_w,j_w},e_j)\times \\
 &&\Big(\psi_i(v'_{i_v,j_v,i_w,j_w}(e_j))+\psi_i(w'_{i_v,j_v,i_w,j_w}(e_j))-
      \psi_i(v_{i_v,j_v})+\psi_i(w_{i_w,j_w})\Big)\,,
\end{eqnarray*}
with
\begin{eqnarray*}
 v'_{i_v,j_v,i_w,j_w}(e_j)&=&
 \frac{v_{i_v,j_v}+w_{i_w,j_w}}{2}+\frac{1}{2}|v_{i_v,j_v}-w_{i_w,j_w}|e_j\,, \\
 w'_{i_v,j_v,i_w,j_w}(e_j)&=&
 \frac{v_{i_v,j_v}+w_{i_w,j_w}}{2}-\frac{1}{2}|v_{i_v,j_v}-w_{i_w,j_w}|e_j\,.
\end{eqnarray*}
It is clear that it impossible to compute all these matrices by the direct use 
of the above formulae for reasonable discretisation parameters. However, the 
separated structure of the factors allows to precompute three arrays 
$P_{GL},P_L$ and $P_Q$ and to use them to assemble the collision matrices in 
an efficient manner.
The components of the first array are for $k=0,\dots,K\,,\ \ell=0,\dots,L$ and 
$i=1,\dots,N_{GL}$
$$
 \Big(P_{GL}\Big)_{k,\ell,i}=\omega_i^{GL}
 x_i^{\ell}\mbox{e}^{x_i/2}\,L_{k}^{(\ell+1/2)}(x_i)
$$
leading to $(K+1)(L+1)N_{GL}$ words of computer memory.
The components of the second array are for $\ell=0,\dots,L\,,\ 
m=-\ell,\dots,\ell$ and 
$j=1,\dots,N_L$
$$
 \Big(P_L\Big)_{l,m,j}=\omega_{j}^LY_{\ell,m}(e_j)
$$
leading to $(L+1)^2N_L$ words of computer memory. Finally, the most 
complicated array is
$$
 \Big(P_Q\Big)_{i_v,j_v,i_w,j_w,i}=\big(q_L\big)_i(v_{i_v,j_v},w_{i_w,j_w})
$$
leading to $n\,N_{GL}^2N_L^2$ words of computer memory. The numerical 
cost of the first two arrays is sub-linear in $n$ and requires a negligible 
computer time.
The computation of the third array, however, is more demanding. 
Formally, 
it requires only a linear 
amount of operations with respect to the number of unknowns and quadratic 
with respect to the 
number of integration points. However, these numbers are not independent and in 
order to keep the spectral accuracy, an increase of the number 
of integration points is unavoidable with increasing $n$.
The computations of all $Q_i\,,\ i=1,\dots,n$ is as follows
$$
 \Big(Q_{N_{GL},N_L}\Big)_i[k,\ell]=(2\,\pi)^2
 \sum_{i_v=1}^{N_{GL}}\sum_{j_v=1}^{N_L}\sum_{i_w=1}^{N_{GL}}\sum_{j_w=1}^{N_L}
 \alpha_{i_v,j_v,i_w,j_w}
 \Big(P_Q\Big)_{i_v,j_v,i_w,j_w,i}
$$
where
$$
 \alpha_{i_v,j_v,i_w,j_w}=
 \Big(P_{GL}\Big)_{k_v,l_v,i_v}\Big(P_L\Big)_{l_v,m_v,j_v}
 \Big(P_{GL}\Big)_{k_w,l_w,i_w}\Big(P_L\Big)_{l_w,m_w,j_w}\,.
$$
Note that $\alpha_{i_v,j_v,i_w,j_w}$ is independent of $i$ and, therefore, 
once computed, can be used to update the $i_v,j_v,i_w,j_w$ sum for all 
matrices $i$. This leads to only a few multiplications and additions 
for the 
entries of the collision matrices without evaluation of special functions.
The numerical work and memory, however, is still of the order ${\cal{O}}(n^3)$.
Furthermore, the entries $k,\ell$ of the collision matrices are independent 
from each other and, therefore, their computation can be done in parallel 
by the use of the open MP software without any additional programming effort.
\subsubsection*{Time discretisation}
Once the mass matrix $M$ and all collision matrices $Q_i$ are computed (from 
now on the subscripts $N_{GL}$ and $N_L$ are omitted), we can start a numerical 
solution of the problem. First of all, the initial right hand side 
$\underline{f}_0$ has to be computed
$$
 \big(\underline{f}_0\big)_i=<f_0,\psi_i>\,,\quad i=1,\dots,n\,.
$$
We use the numerical quadrature
$$
 \big(\underline{f}_0\big)_i=2\,\pi\,\sum_{i_v=1}^{N_{GL}}\omega_{i_v}^{GL}
 \mbox{e}^{x_{i_v}}
 \sum_{j_v=1}^{N_L}\omega_{j_v}^L f_0(\sqrt{x_{i_v}}\,e_{j_v})\,
 \psi_i(\sqrt{x_{i_v}}\,e_{j_v})\,,\quad i=1,\dots,n\,
$$
and compute the initial coefficient vector
$$
 \underline{f}^{(0)}=M^{-1}\underline{f}_0\,.
$$
Then we choose a time step $\tau>0$ and the final time $T=\tau\,N_t$.
For the time integration of the system (\ref{Eqn ODE System}), we can 
choose any classical solver, for example the most simple Euler scheme or 
the Runge-Kutta method of the second or a higher order.
For the Euler method the $k$th step, $k=0,\dots,N_t-1$, is as follows.
Compute the vector $\underline{q}$ with
$$
 \underline{q}_i=
 (Q_i\underline{f}^{(k)},\underline{f}^{(k)})\,,\quad 
 i=1,\dots,n
$$
by the use of the BLAS library. Compute the next coefficient vector
$\underline{f}^{(k+1)}$ as
$$
 \underline{f}^{(k+1)}=\underline{f}^{(k)}+\tau\,M^{-1}\underline{q}
$$
by utilising the functionality of the LAPACK package and the BLAS library 
once again.
\section{Numerical examples}
In this section we consider three examples of relaxation. The collision 
kernel of the first two examples will be constant, i.e.
$$
 B(v,w,e)=\dfrac{1}{4\,\pi}\,.
$$
This is the most simple case of Maxwell pseudo--molecules. For this kernel,
the exact relaxation time of any moment of the distribution function is 
known.
Thus, we will be able to check the accuracy of our scheme very carefully. In 
the first example, we will consider a sum of two Maxwell distributions as 
an initial condition. The second example is the famous BKW solution for which 
not only the time relaxation of the moments but the distribution function 
itself is analytically known. However, this solution is an isotropic
function and, therefore, its numerical approximation by our spectral scheme is 
rather simple. The third example will be the classical hard spheres model with 
the collision kernel
$$
 B(v,w,e)=\dfrac{1}{4\,\pi}\,|v-w|\,.
$$
No analytic dependencies are available for this example.
Thus, we will compare our 
results with those obtained by the use of a stochastic particle method. 
\subsection{Relaxation of a mixture of two Maxwellian's}
For the spatially homogeneous relaxation, the density, the mean velocity and 
the temperature
$$ 
 \rho=\int\limits_{\mathbb{R}^3}f(t,v)\,dv\,,\quad
 V=\frac{1}{\rho}\int\limits_{\mathbb{R}^3}v\,f(t,v)\,dv\,,\quad
 T=\frac{1}{3\,\rho}\int\limits_{\mathbb{R}^3}|v-V|^2\,f(t,v)\,dv
$$
are conserved quantities. The relaxation of the flow of momentum, 
the flow of energy and of the special fourth moment
\begin{equation*}
 M(t)=\int\limits_{\mathbb{R}^3}vv^\top\,f(t,v)\,dv\,,\ 
 r(t)=\int\limits_{\mathbb{R}^3}v|v|^2\,f(t,v)\,dv\,,\ 
 s(t)=\int\limits_{\mathbb{R}^3}|v|^4\,f(t,v)\,dv
\end{equation*}
is given as in \cite{RJASANOW;WAGNER:05} by
\begin{eqnarray}
\label{Eqn Exact Moments}
 M(t)&=&M_0\,e^{-t/2}+\Big(T\, I+VV^\top\Big)\Big(1-e^{-t/2}\Big)\,, 
 \\
 r(t)&=&r_0\,e^{-t/3}+\Big(5\,T+|V|^2\Big)\,V\,\Big(1-e^{-t/3}\Big) \nonumber \\
 &&+
 2\Big(M_0-VV^\top-T\, I\Big)\,V\,\Big(e^{-t/2}-e^{-t/3}\Big)\,, \nonumber \\
 \label{Eqn Exact Moments 2}
 s(t)&=&s_0\,e^{-t/3}+\Big(|V|^4+15\,T^2+10\,T\,|V|^2\Big)\Big(1-e^{-t/3}\Big)  
 \\
 &&+\frac{1}{2}\Big(||M_0||_F^2-3\,T^2+|V|^4-2\big( M_0V,V \big)\Big)
 \Big(e^{-t}-e^{-t/3}\Big)  \nonumber \\
 &&+4\,\Big(\big( M_0V,V \big)-|V|^4-T\,|V|^2\Big)
 \Big(e^{-t/2}-e^{-t/3}\Big)\,, \nonumber
\end{eqnarray}
where
\begin{equation*}
M_0=\int\limits_{\mathbb{R}^3}vv^\top\,f_0(v)\,dv\,,\quad
r_0=\int\limits_{\mathbb{R}^3}v|v|^2\,f_0(v)\,dv\,,\quad
s_0=\int\limits_{\mathbb{R}^3}|v|^4\,f_0(v)\,dv
\end{equation*}
and $\|\cdot\|_F$ denotes the Frobenius norm.
We will consider the initial distribution $f_0$ in the form of a convex sum of 
two Maxwell distributions
\begin{equation}
 f_0(v)=\alpha f_{M_1}(v)+(1-\alpha)f_{M_2}(v)\,,\quad 0\le\alpha\le 1\,,
\end{equation}
where
$$
 f_{M_i}(v)=\frac{1}{(2\pi\,T_i)^{3/2}}\,
 \mbox{e}^{\displaystyle-\frac{|v-V_i|^2}{2\,T_i}}\,,\quad i=1,2\,.
$$
In these settings, the initial values are
\begin{eqnarray*}
 \rho&=&1\,, \\
 V&=&\alpha V_1 + (1-\alpha ) V_2\,, \\
 T&=&\alpha T_1+(1-\alpha )T_2+
 \frac{1}{3}\alpha (1-\alpha )|V_1-V_2|^2\,, \\
 M_0&=&\alpha\Big(T_1\,I+V_1V_1^\top\Big)+
 (1-\alpha)\Big(T_2\,I+V_2V_2^\top\Big)\,, \\
 r_0&=&\alpha\Big(5T_1+|V_1|^2\Big)V_1+
 (1-\alpha )\Big(5T_2+|V_2|^2\Big)V_2\,, \\
 s_0&=&\alpha\Big(|V_1|^4+15\,T_1^2+10\,T_1\,|V_1|^2\Big)+ \\ 
 &&
 (1-\alpha )\Big(|V_2|^4+15\,T_2^2+10\,T_2\,|V_2|^2\Big)\,.
\end{eqnarray*}
For our first example, we choose
$$
 \alpha=1/2\,,\ V_1=(-1,0,0)^\top\,,\quad V_2=(+1,0,0)^\top\,,\quad 
T_1=T_2=\frac{2}{3} 
$$
and obtain
$$
 V=(0,0,0)^\top\,,\quad T=1\,.
$$
\newpage
\subsubsection*{Initial condition}
\begin{figure}
\centering
\includegraphics[width=13cm]{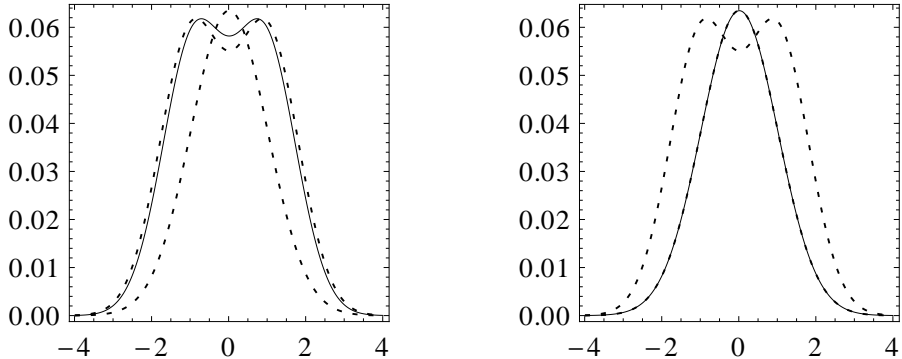}
\caption{Initial and final distributions for $n=27$}
\label{Fig Initial and final n=27}
\end{figure}

\begin{figure}
\centering
\includegraphics[width=13cm]{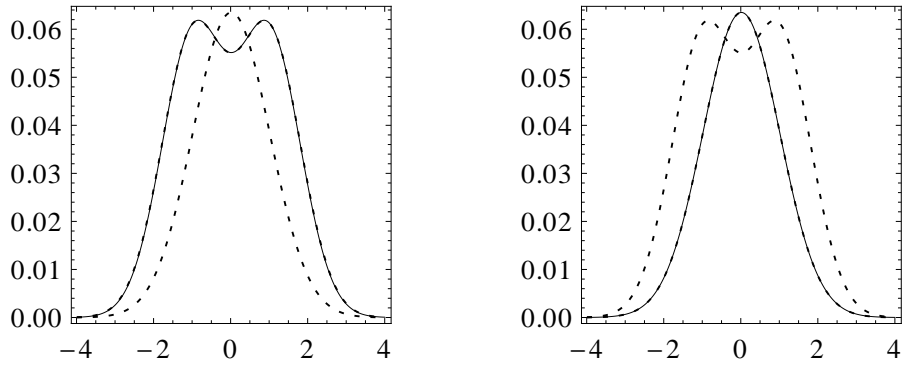}
\caption{Initial and final distributions for $n=125$}
\label{Fig Initial and final n=125}
\end{figure}

\begin{figure}
\centering
\includegraphics[width=13cm]{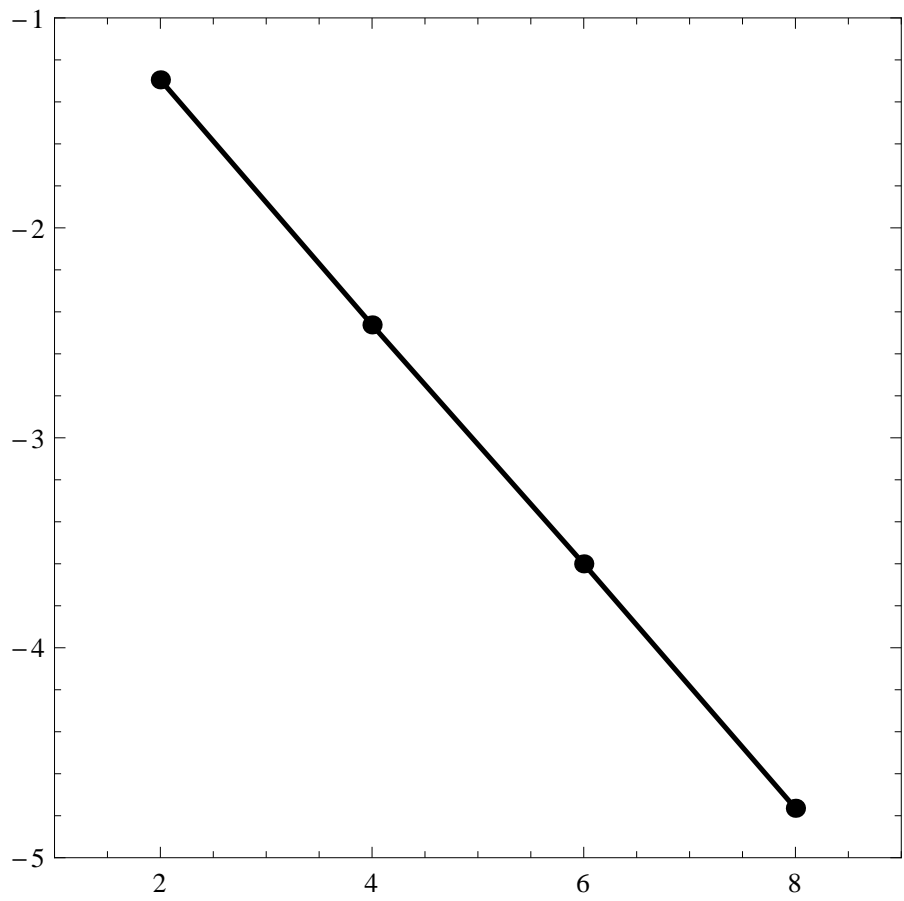}
\caption{$\log_{10}$ course of the $L_2(\mathbb{R}^3)$ error}
\label{Fig Log L2 error}
\end{figure}

For a series of discretisation parameters $K$ and $L$, we first define the 
parameters of the Gauss-Laguerre and Lebedev quadratures $N_{GL}$ and $N_L$
in the following way.
We perform approximation of the initial condition and choose the minimal values 
of $N_{GL}$ and $N_L$ leading to the 
highest approximation quality for the given values of $K$ and $L$. 
In the first two figures we illustrate the approximation of the initial 
condition $f_0(v_1,0,0)\,,\ v_1\in[-4,4]$ and of the final Maxwell distribution 
$f_M(v_1,0,0)\,,\ v_1\in[-4,4]$ for $K=L=2$ with 
$n=27$ basis functions (Figure \ref{Fig Initial and final n=27}) and for $K=L=4$ 
with 
$n=125$ basis functions (Figure \ref{Fig Initial and final n=125}). The initial 
condition and the final Maxwell distribution are shown with thick dashed lines, 
while the numerical approximation is depicted by the thin solid line. There 
is 
a clear numerical error by the approximation of the initial condition for  
$n=27$. For $n=125$, however, the error can not be optically seen on the 
figure. The final Maxwell distribution is perfectly approximated in both cases.
The $L_2(\mathbb{R}^3)$ error 
$$
\dfrac{\|f^{(n)}(0,\cdot)-f_0\|_{L_2(\mathbb{R}^3)}}
{\|f_0\|_{L_2(\mathbb{R}^3)}}
$$
of the approximation of the initial condition 
$f_0$ is summarised in Table \ref{Tbl Initial error} and its logarithmic plot is 
shown in Figure \ref{Fig Log L2 error}. The last column in Table \ref{Tbl 
Initial error} contains the Convergence Factor (CF), i.e. a quotient of two 
consecutive errors. The exponential convergence of the error 
is clearly seen.
\begin{table}
\centering
\caption{Approximation error for the initial condition}
\begin{tabular*}{\textwidth}{@{\extracolsep{\fill}}rrrrrrr}
\hline\noalign{\smallskip}
 $K$ & $L$ & $n$ & $N_L$ & $N_{GL}$ & $L_2(\mathbb{R}^3)$-Norm & CF \\
\noalign{\smallskip}\hline\noalign{\smallskip}
 2 & 2 &   27 &  38 &  8 & $5.07\cdot10^{-2}$ &       - \\ 
 4 & 4 &  125 &  50 &  8 & $3.45\cdot10^{-3}$ &    14.7 \\ 
 6 & 6 &  343 & 110 & 16 & $2.51\cdot10^{-4}$ &    13.7 \\ 
 8 & 8 &  729 & 110 & 16 & $1.72\cdot10^{-5}$ &    14.6 \\ 
\noalign{\smallskip}\hline
\end{tabular*}
\label{Tbl Initial error}
\end{table}
\newpage
\subsubsection*{Relaxation of the moments}
\begin{figure}
\centering
\includegraphics[width=13cm]{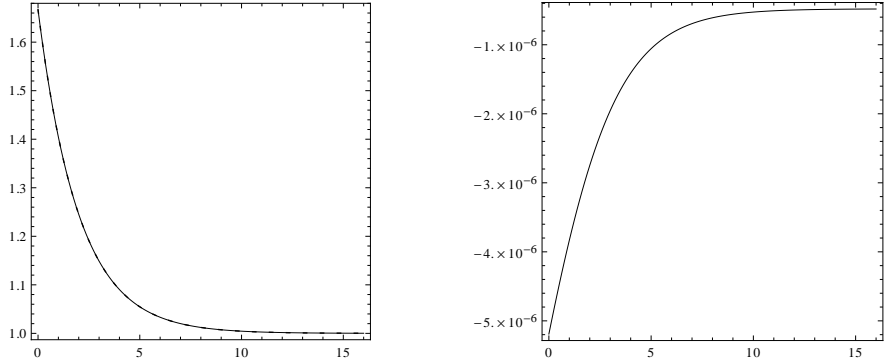}
\caption{Course of the functional $M_{11}(t)$}
\label{Fig M11}
\end{figure}

\begin{figure}
\centering
\includegraphics[width=13cm]{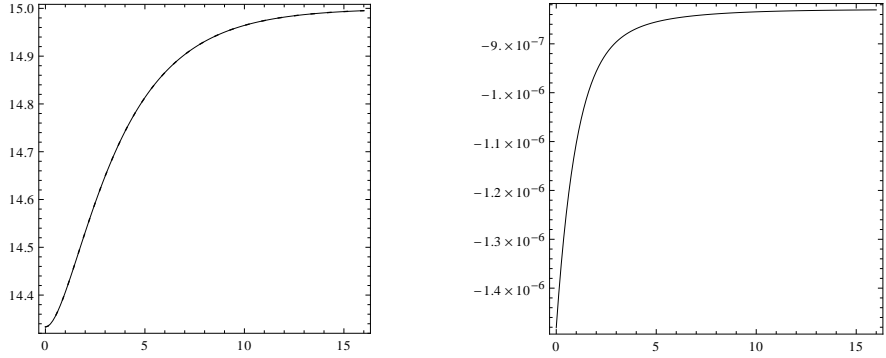}
\caption{Course of the functional $s(t)$}
\label{Fig s}
\end{figure}

In the study of the accuracy of the time dependent moments of the 
distribution function, two new aspects have to be considered, namely the value 
of the time 
discretisation parameter $\tau$ and the quality of the time integrating scheme. 
We will demonstrate the efficiency of the simplest Euler scheme and of the 
Runge-Kutta method of orders two and four. For the given example, there is 
a non-trivial relaxation of the main diagonal components of the flux of 
momentum 
tensor (\ref{Eqn Exact Moments}) and of the fourth moment (\ref{Eqn Exact 
Moments 2}). Figure \ref{Fig M11} shows the course of the function $M_{11}(t)$ 
where the thick dashed line is the analytic solution and the thin solid line 
is the computed moment for $n=125$. There is no optical difference. The time 
relaxation of the function $s(t)$ is shown in Figure \ref{Fig s}. The right 
plots on both figures show the time evolution of the difference between the 
analytic and the numerical solutions.
In the next three tables we 
present the maximal error for these moments computed with different time steps 
on the time interval $[0,16]$ for different values of $n$. The lines indicated 
in bold, shows the best accuracy reached for the given value of $n$ and for a 
maximal number of time steps $N_t\le 8192$. The linear, 
quadratic and fourth order convergences in time for different time integration 
schemes are evident.
\begin{table}
\centering
\caption{Error for the moments $M_{11}(t)$ and $s(t)$, Euler method}
\begin{tabular*}{\textwidth}{@{\extracolsep{\fill}}cccccc}
\hline\noalign{\smallskip}
 $N_t$ & $n$ & $M_{11}(t)$ & CF & $s(t)$ & CF \\
\noalign{\smallskip}\hline\noalign{\smallskip}
  32 &  27 & $2.78\cdot10^{-2}$ &    - & $1.75\cdot10^{-3}$ &    - \\ 
  64 &  27 & $1.33\cdot10^{-2}$ & 2.09 & $8.01\cdot10^{-4}$ & 2.18 \\ 
 128 &  27 & $6.61\cdot10^{-3}$ & 2.01 & $4.13\cdot10^{-4}$ & 1.94 \\ 
 {\bf 256} & {\bf 27} & ${\bf 3.37\cdot10^{-3}}$ & {\bf 1.96} & 
                        ${\bf 2.40\cdot10^{-4}}$ & {\bf 1.72} \\ 
 256 & 125 & $3.17\cdot10^{-3}$ &    - & $1.61\cdot10^{-4}$ &    - \\ 
 512 & 125 & $1.58\cdot10^{-3}$ & 2.01 & $7.98\cdot10^{-5}$ & 2.02 \\ 
1024 & 125 & $7.89\cdot10^{-4}$ & 2.00 & $4.01\cdot10^{-5}$ & 1.99 \\ 
2048 & 125 & $3.95\cdot10^{-4}$ & 2.00 & $2.05\cdot10^{-5}$ & 1.96 \\ 
4096 & 125 & $1.99\cdot10^{-4}$ & 1.98 & $1.08\cdot10^{-5}$ & 1.90 \\ 
{\bf 8192} & {\bf 125} & ${\bf 1.00\cdot10^{-4}}$ & {\bf 1.99} & 
                         ${\bf 5.98\cdot10^{-6}}$ & {\bf 1.81} \\ 
\noalign{\smallskip}\hline
\end{tabular*}
\label{Tbl Moments Error RK Euler}
\end{table}
\begin{table}
\centering
\caption{Error for the moments $M_{11}(t)$ and $s(t)$, Runge-Kutta 2}
\begin{tabular*}{\textwidth}{@{\extracolsep{\fill}}cccccc}
\hline\noalign{\smallskip}
 $N_t$ & $n$ & $M_{11}(t)$ & CF & $s(t)$ & CF \\
\noalign{\smallskip}\hline\noalign{\smallskip}
   32 &  27 & $2.34\cdot10^{-3}$ &    - & $2.63\cdot10^{-4}$ &    - \\ 
  {\bf 64} & {\bf 27} & ${\bf 4.16\cdot10^{-4}}$ & {\bf 5.63} & 
                       ${\bf 9.93\cdot10^{-5}}$ & {\bf 4.19} \\ 
   64 & 125 & $5.72\cdot10^{-4}$ &    - & $7.49\cdot10^{-5}$ &    - \\ 
  128 & 125 & $1.34\cdot10^{-4}$ & 4.27 & $1.68\cdot10^{-5}$ & 4.46 \\ 
  256 & 125 & $3.11\cdot10^{-5}$ & 4.31 & $3.26\cdot10^{-6}$ & 5.15 \\ 
  {\bf 512} & {\bf 125} & ${\bf 6.03\cdot10^{-6}}$ & {\bf 5.16} & 
                         ${\bf 1.48\cdot10^{-6}}$ & {\bf 2.20} \\ 
 1024 & 343 & $2.05\cdot10^{-6}$ &    - & $2.63\cdot10^{-7}$ &    - \\ 
 2048 & 343 & $5.12\cdot10^{-7}$ & 4.00 & $6.55\cdot10^{-8}$ & 4.05 \\ 
 4096 & 343 & $1.27\cdot10^{-7}$ & 4.03 & $1.63\cdot10^{-8}$ & 4.02 \\ 
  {\bf  8192} & {\bf 343} & {\bf $3.19\cdot10^{-8}$} & {\bf 3.98} &
                           {\bf $4.07\cdot10^{-9}$} & {\bf 4.00} \\ 
\noalign{\smallskip}\hline
\end{tabular*}
\label{Tbl Moments Error RK 2}
\end{table}
\begin{table}
\centering
\caption{Error for the moments $M_{11}(t)$ and $s(t)$, Runge-Kutta 4}
\begin{tabular*}{\textwidth}{@{\extracolsep{\fill}}cccccc}
\hline\noalign{\smallskip}
 $N_t$ & $n$ & $M_{11}(t)$ & CF & $s(t)$ & CF \\
\noalign{\smallskip}\hline\noalign{\smallskip}
  {\bf 32} & {\bf 27} & ${\bf 4.16\cdot10^{-4}}$ & {\bf -} & 
                       ${\bf 4.20\cdot10^{-4}}$ & {\bf -} \\ 
  {\bf 32} & {\bf 125} & ${\bf 5.81\cdot10^{-6}}$ & {\bf -} & 
                         ${\bf 5.41\cdot10^{-6}}$ & {\bf -} \\ 
   64 & 343 & $4.53\cdot10^{-7}$ &    - & $2.18\cdot10^{-8}$ &    - \\ 
  128 & 343 & $2.68\cdot10^{-8}$ & 16.9 & $1.28\cdot10^{-9}$ & 17.0 \\ 
  256 & 343 & $1.58\cdot10^{-9}$ & 17.0 & $6.33\cdot10^{-11}$ & 20.2 \\ 
  {\bf 512} & {\bf 343} & ${\bf 1.36\cdot10^{-10}}$ & {\bf 9.41} & 
                         ${\bf 2.43\cdot10^{-11}}$ & {\bf 2.61} \\ 
  128 & 729 & $2.68\cdot10^{-8}$ & - & $1.28\cdot10^{-9}$ &  - \\ 
  256 & 729 & $1.58\cdot10^{-9}$ & 17.0 & $6.32\cdot10^{-11}$ & 20.2 \\ 
  {\bf 512} & {\bf 729} & ${\bf 1.71\cdot10^{-14}}$ & {\bf -} & 
                         ${\bf 1.63\cdot10^{-9}}$ & {\bf  -} \\ 
\noalign{\smallskip}\hline
\end{tabular*}
\label{Tbl Moments Error RK 4}
\end{table}
In the last table, we observe no proper convergence of the finest 
discretisation with $n=729$. The errors for $N_t=128,256$ is 
practically identical to those obtained for $n=343$. For $N_t=512$, the error 
for $M_{11}$ practically jumps to the machine accuracy, while the error for 
$s$ increases.
This is a clear indicator that the numerical integration with $N_{LG}=16$ 
and $N_L=110$ is not sufficiently accurate to yield the theoretically 
achievable high accuracy for this $n$.
\subsubsection*{$H$-functional and convergence to equilibrium}
\begin{figure}
\centering
\includegraphics[width=13cm]{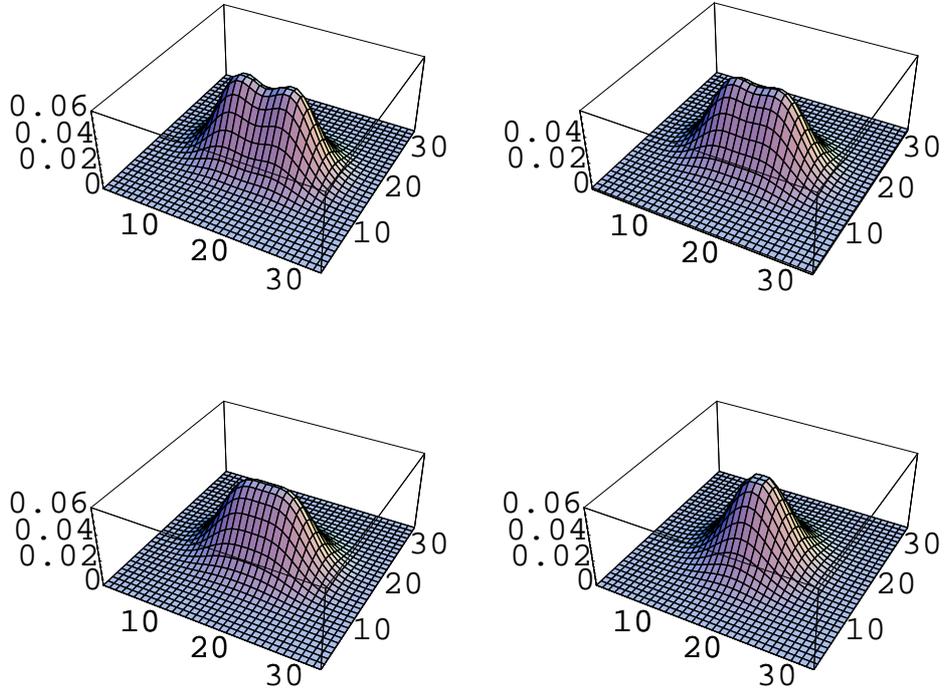}
\caption{Density function for $t=0,1/4,1,16$}
\label{Fig Density plots}
\end{figure}

\begin{figure}
\centering
\includegraphics[width=13cm]{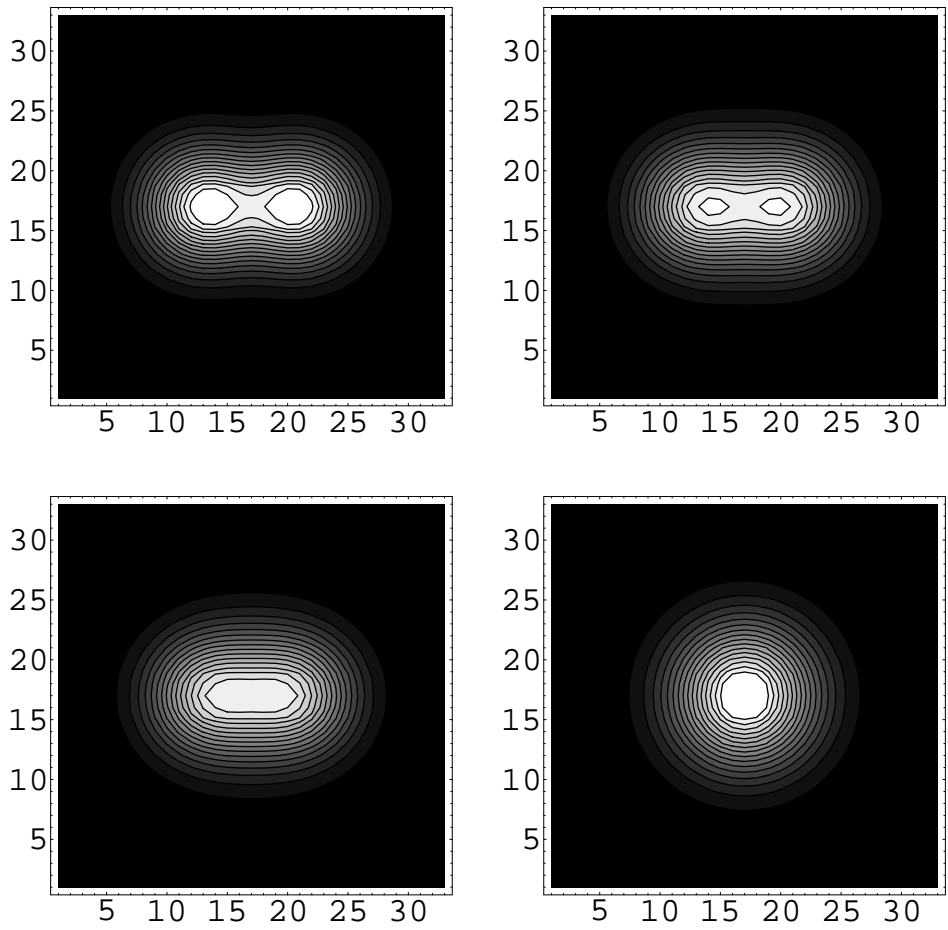}
\caption{Contours of the density function for $t=0,1/4,1,16$}
\label{Fig Conrours plots}
\end{figure}

\begin{figure}
\centering
\includegraphics[width=13cm]{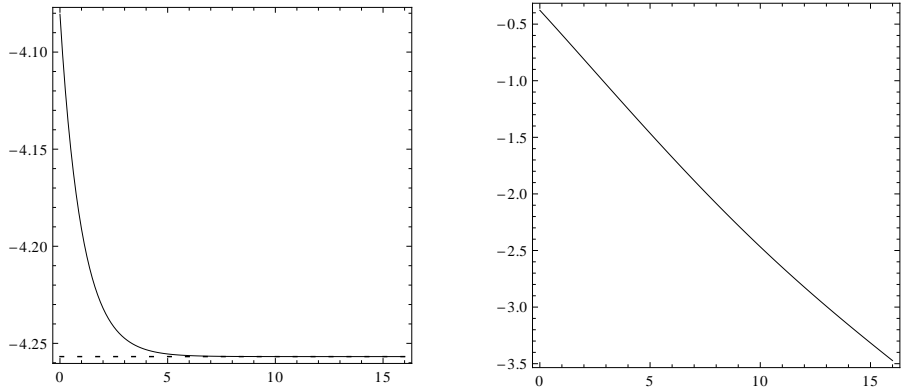}
\caption{Course of the $H$-functional and of the $L_2-$error for $n=343$}
\label{Fig H Functional}
\end{figure}

In Figures \ref{Fig Density plots} and \ref{Fig Conrours plots},
we show the plots of the numerical density function $f^{(n)}(t,v_1,v_2,0)$
and its contours
for $(v_1,v_2)\in[-4,4]\times[-4,4]$ with $32\times 32$ points and for 
the times $t=0,1/4,1,16$ obtained for $n=125$.
\newpage
Finally, in Figure \ref{Fig H Functional}, we show the time relaxation of the 
numerical Boltzmann $H$-functional (left plot)
$$
 H(t)=\int\limits_{\mathbb{R}^3}f^{(n)}(t,v)\ln f^{(n)}(t,v)\,dv
$$
for $n=125$. Its analytically known asymptotic value
$$
 \lim_{t\rightarrow\infty}H(t)=
 \ln\dfrac{1}{(2\,\pi)^{3/2}}-\dfrac{3}{2}=-4.25681\dots\,.
$$
is shown as a dashed thick line, while the course of the $H$-functional is 
depicted by the 
thin solid line. The right plot in Figure \ref{Fig H Functional} shows the 
$\log_{10}$-course of the relative $L_2$-norm of the difference of the 
current distribution function to the final Maxwell distribution.
$$
 \dfrac{\|f^{(n)}(t,\cdot)-f_M\|_{L_2(\mathbb{R}^3)}}
{\|f_M\|_{L_2(\mathbb{R}^3)}}
$$
which obviously shows exponential convergence.
\subsection{BKW solution}
In this subsection, we consider the famous exact solution of the Boltzmann 
equation found by Bobylev \cite{BOBYLEV:75} and Krook and Wu \cite{KROOK;WU:77}.
The solution is obtained for $\lambda=0$ in (\ref{Eqn Collision Kernel}) and 
is of the form
\begin{eqnarray}
\lefteqn{
 f(t,v)= }  \nonumber \\
&&
\frac{\rho}{(2\,\pi T)^{3/2}}(\beta(t)+1)^{3/2}
 \left(1+\beta(t)\Big(\frac{\beta(t)+1}{2T}|v|^2-\frac{3}{2}\Big)\right)
 \mbox{e}^{\displaystyle -\frac{\beta(t)+1}{2T}|v|^2}\,, \nonumber
\end{eqnarray}
with
$$
 \beta(t)=\frac{\beta_0\,\mbox{e}^{\displaystyle -\alpha\,\rho\, t/2}}
 {1+\beta_0\,(1-\mbox{e}^{\displaystyle-\alpha\,\rho\, t/2})}\,,
$$
where $\beta_0$ denotes the initial value for the function 
$\beta$ and $\alpha$ is defined as
$$
 \alpha=C_0\pi\int\limits_0^{\pi}b(\cos\theta)\sin^3\theta\,d\theta\,.
$$
This solution is non-negative for 
$$
0 \le \beta_0 \le 2/3\,.
$$
The density $\rho$ and the temperature $T$ are two additional parameters.
We will use the following setting for our tests
$$
 C_0=\dfrac{1}{4\,\pi}\,,\ b(\cos\theta)=1\,,\ 
 \alpha=1/3\,,\ \rho=1\,,\ T=1\,,\ \beta_0=2/3, 
$$
leading to the solution
$$
 f(t,v)\!=\! 
\frac{1}{(2\,\pi)^{3/2}}\,\big(\beta(t)+1\big)^{3/2}\,
 \left(1+\beta(t)\Big(\frac{\beta(t)+1}{2}|v|^2-\frac{3}{2}\Big)\right)
 \mbox{e}^{\displaystyle -\frac{\beta(t)+1}{2}|v|^2}\,,
$$
where
$$
 \beta(t)=\frac{2\,\mbox{e}^{\displaystyle -t/6}}
 {5-2\,\mbox{e}^{\displaystyle -t/6}}\,.
$$
\subsubsection*{Initial condition}
\begin{figure}
\centering
\includegraphics[width=13cm]{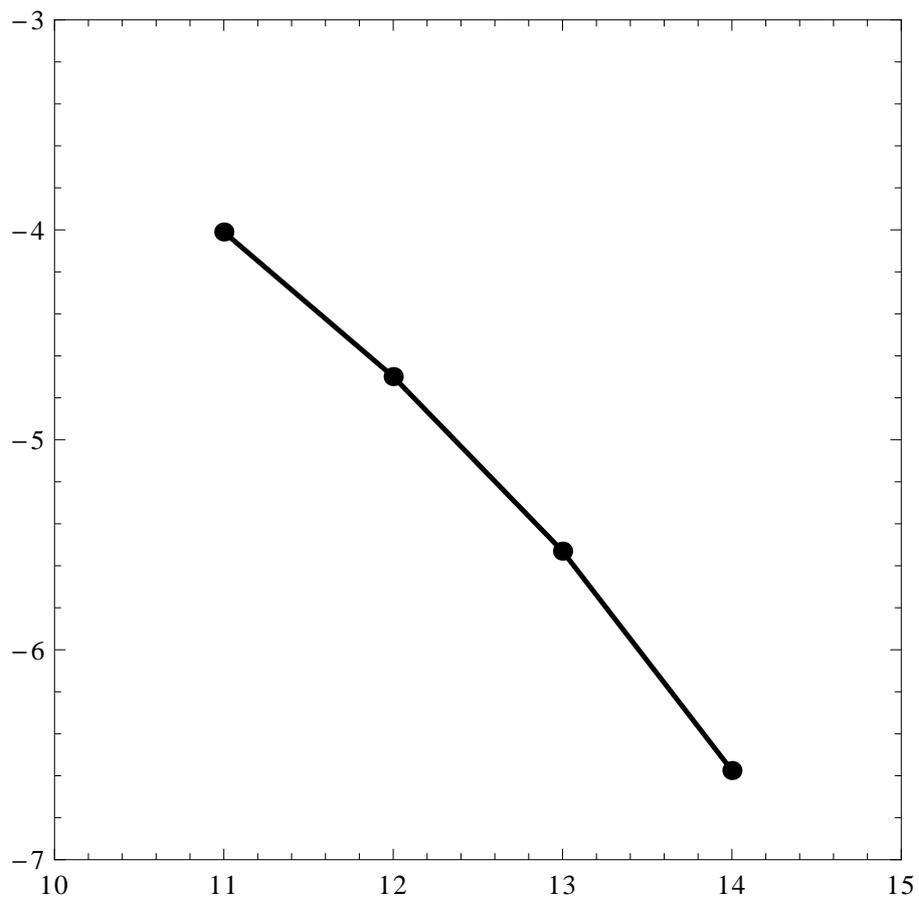}
\caption{$\log_{10}$ course of the $L_2(\mathbb{R}^3)$ error, BKW}
\label{Fig Log L2 error BKW}
\end{figure}

Since the BKW solution is an isotropic function, we change only the parameter 
$K$
and let $L$ be zero for all tests. This leads to a very low number of unknowns 
and to an extremely fast numerical solution of the Boltzmann equation 
taking only a few seconds.
A stable spectral convergence starts with $K=11$ and the results of the 
approximation of the initial condition are shown in Table \ref{Tbl Initial 
error BKW} and in Figure \ref{Fig Log L2 error BKW}. However, for these values 
no optical difference to the initial condition and to the final Maxwell 
distribution can be seen on a figure. Thus, we show the approximation of the 
initial condition and of the final Maxwell distribution for $K=6\,,\ n=7$ in 
Figure \ref{Fig Initial and final n=7, BKW}.
\begin{figure}
\centering
\includegraphics[width=13cm]{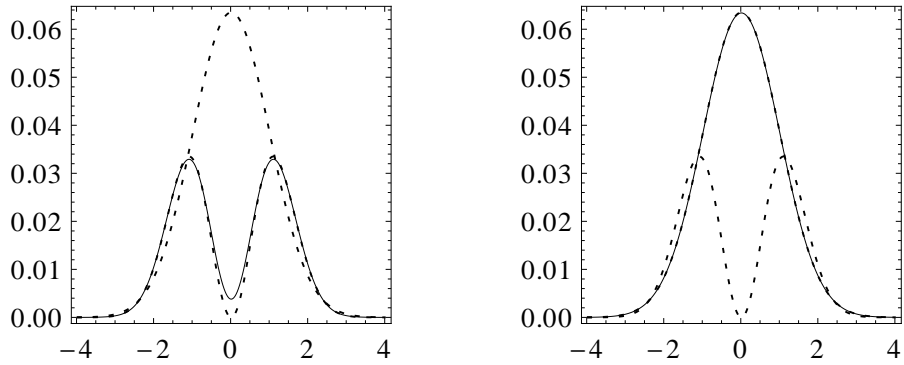}
\caption{Initial and final distributions for $n=7$, BKW}
\label{Fig Initial and final n=7, BKW}
\end{figure}
\begin{table}
\centering
\caption{Approximation error for the BKW initial condition}
\begin{tabular*}{\textwidth}{@{\extracolsep{\fill}}rrrrrrr}
\hline\noalign{\smallskip}
 $K$ & $L$ & $n$ & $N_L$ & $N_{GL}$ & $L_2(\mathbb{R}^3)$-Norm & CF \\
\noalign{\smallskip}\hline\noalign{\smallskip}
 11 & 0 & 12 & 38 & 16 & $9.77\cdot10^{-5}$ &     - \\ 
 12 & 0 & 13 & 38 & 16 & $2.00\cdot10^{-5}$ &  4.89 \\ 
 13 & 0 & 14 & 38 & 16 & $2.95\cdot10^{-6}$ &  6.78 \\ 
 14 & 0 & 15 & 38 & 16 & $2.66\cdot10^{-7}$ & 11.09 \\ 
\noalign{\smallskip}\hline
\end{tabular*}
\label{Tbl Initial error BKW}
\end{table}
\subsubsection*{Relaxation of the moments}
For the BKW solution, all physical moments remain constant in time and they 
are approximated with an accuracy of about $10^{-14}-10^{-15}$ even for $n=7$.
Thus, we show only the course of the fourth moment $s(t)$ on the time interval 
$[0,16]$ as well as the difference to the exact curve for $n=7$ in Figure 
\ref{Fig s, BKW}. The results are obtained with the Runge-Kutta method of the 
fourth order with $N_t=128$ time steps.
\begin{figure}
\centering
\includegraphics[width=13cm]{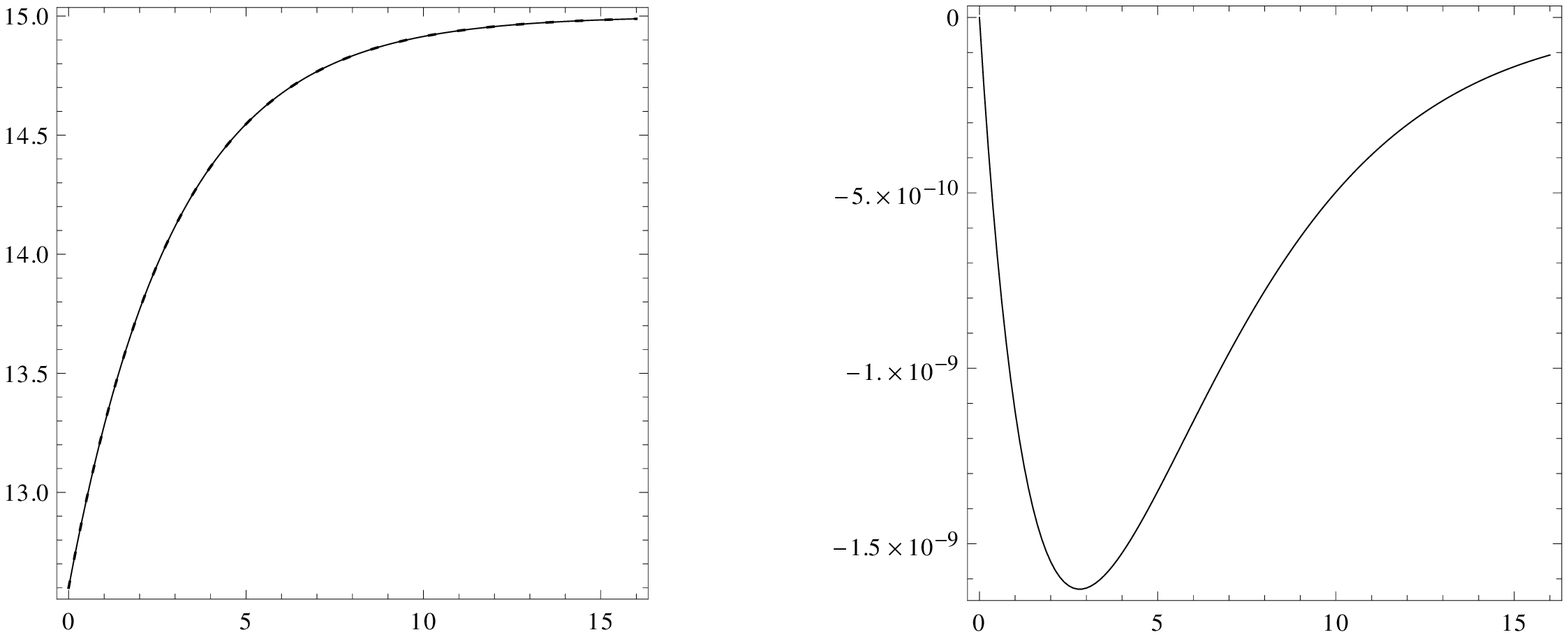}
\caption{Course of the functional $s(t)$, BKW}
\label{Fig s, BKW}
\end{figure}
\subsection{Hard spheres}
There is no analytic information about the exact solution for the case of 
hard spheres. Thus, we 
consider the above mixture of two Maxwell distributions as the initial 
condition and choose the solution obtained by the use of the stochastic 
particle method (see \cite{RJASANOW;WAGNER:05}) as a reference. We 
choose $8192$ equally weighted particles and compute $8192$ independent 
trajectories of the process on the time interval $[0,4]$. Thus the accuracy 
of the stochastic solution should be of the order $10^{-3}-10^{-4}$.
For comparison, we 
take the curves obtained for $K=L=4$, i.e. for $n=125$ unknowns. For the time 
integration, the Runge-Kutta method of the fourth order with $N_t=128$ time 
steps has been used. The dependence of the moment $M_{11}(t)$ on time
is shown in Figure 
\ref{Fig M11, HS}. The thick dashed line represents the stochastic reference 
solution on the left plot. The thin solid line is the Galerkin-Petrov solution. 
The right plot shows the difference between the curves. The accuracy is of 
the order $10^{-4}$. The same data is shown for the fourth moment $s(t)$ 
in Figure \ref{Fig s, HS}. Here, the accuracy is of the order $10^{-3}$ and 
some oscillations of the stochastic solution are apparent. The 
computational time for the stochastic particle method
on a single Intel i7 processor was about 10 minutes while the Galerkin-Petrov 
solution with precomputed collision matrices ($N_{GL}=8,\ N_L=50$) was obtained 
in 20 seconds. The computational time for the collision matrices was about 12 
minutes, which is to the computation time of the stochastic solution. 
It seems, that the error is mostly due to the stochastic 
approximation, but to obtain an additional order of its accuracy, 
the computational effort must be increased 100-fold .
\begin{figure}
\centering
\includegraphics[width=13cm]{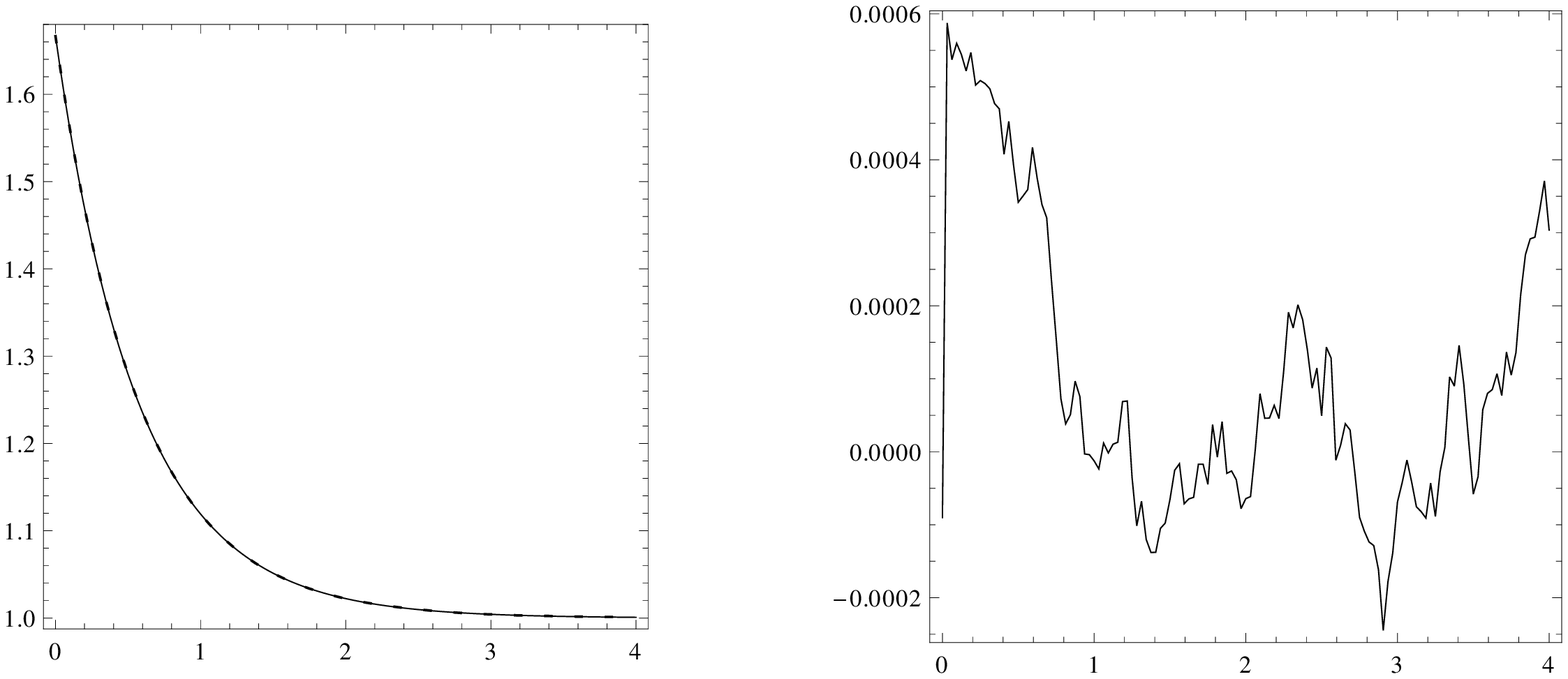}
\caption{Course of the functional $M_{11}(t)$, Hard Spheres}
\label{Fig M11, HS}
\end{figure}
\begin{figure}
\centering
\includegraphics[width=13cm]{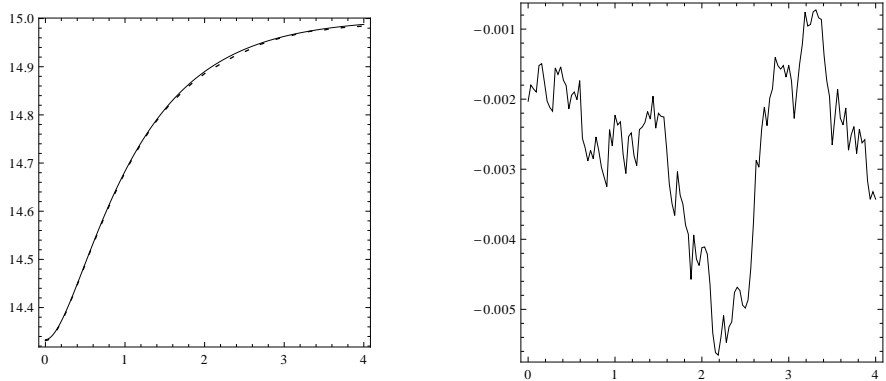}
\caption{Course of the functional $s(t)$, Hard Spheres}
\label{Fig s, HS}
\end{figure}
\section{Conclusions}
In this paper we present a new deterministic numerical scheme for the classical 
spatially homogeneous Boltzmann equation. The scheme is based on a spectral 
Galerkin-Petrov scheme. The main features of 
the method are the following:
\begin{enumerate}
 \item The method uses mutually orthonormal, globally defined basis functions 
       derived from the normalised Maxwell distribution, the Laguerre 
       polynomials and the spherical harmonics;
 \item The system of test functions consists of globally defined low order 
       polynomials;
 \item Since the set of test functions contains all collision invariants, the 
       method is automatically conservative, i.e. the numerical collision 
       invariants remains constant up to the machine accuracy without any 
       additional numerical effort;
 \item The approximation quality of the method is spectral, i.e. there is an 
       exponential convergence. However, this property holds only for 
       infinitely smooth functions;
 \item The main numerical work of the method is the initial computation of 
       the collision matrices. However, once computed, these matrices can 
       be used for different initial conditions, on different time 
       intervals and for different time integration schemes. Then the 
       computational procedure for the whole relaxation in time takes only 
       seconds on a single processor;
 \item Two classical numerical examples for the spatially homogeneous 
       relaxation, namely mixture of two Maxwell distributions as an 
       initial condition and the BKW solution were computed up to a very high 
       accuracy with a low number of basis function of $10^1-10^3$;
 \item The error due to the time integration dominated over the
       spectral error. The Runge-Kutta method of the fourth order was 
       sufficient to equalise both errors;
 \item For the hard spheres model, we've shown an excellent agreement of the 
       results obtained by the new scheme with those obtained by a stochastic 
       particle scheme.
\end{enumerate}
A future work in this research area should certainly contain the following 
points:
\begin{enumerate}
 \item Development of numerical integration quadratures for non-cutoff kernels 
       $B$ for an effective evaluation of the integrals 
       (\ref{Eqn Collision Distribution}), as mentioned in the end of 
       Section \ref{Sec Galerkin-Petrov approximation};
 \item Spatially homogeneous numerical tests to understand how far can 
       the computations be done with the deviation of the temperature from 
       its value equal to one and with the deviation from the zero mean 
       velocity. This will help to formulate criteria for an enrichment 
       of the system of the basis functions;
 \item The method can be easily adapted for the inelastic Boltzmann equation 
       with constant or even variable (relative velocity dependent) restitution 
       coefficient. The main difference is in the term (\ref{Eqn Collision 
       Distribution}). However, since the tails of the distribution function
       of the inelastic Boltzmann equation exhibit an asymptotic different 
       from the Maxwell distribution, the system of basis function should be 
       modified as well;
 \item The main goal is an application of the proposed approach to the 
       spatially inhomogeneous Boltzmann equation. In this case the system 
       (\ref{Eqn ODE System}) of ODE's will be transformed into the 
       hyperbolic system
\begin{equation}
 \label{Eqn Hyperbolic System}
 \dfrac{\partial}{\partial t}\Big(M\underline{f}(t,x)\Big)_i+
 \mbox{div}_x(F_i\underline{f}(t,x))=
 \underline{f}(t,x)^\top\,Q_i\,\underline{f}(t,x)\,,\quad
 i=1,\dots,n\,, \nonumber
\end{equation}
       where the flow matrices $F_i\in\mathbb{R}^{3\times n}$
       have the entries
$$
 F_i[m,j]=<v_m\varphi_j,\psi_i>\,,\quad m=1,2,3\,,\ j=1,\dots,n\,.
$$
       for $i=1,\dots,n$. These matrices can be easily precomputed and stored 
       requiring much less memory than the collision matrices $Q_i$;
\item The proposed spatially inhomogeneous method can be especially efficient 
      for very slow flows with a small deviation of the temperature from its 
      mean value. Exactly for such flows, the application of the stochastic 
      particle methods is problematic;
\item In spatially inhomogeneous flows, the situations occur where 
      the distribution function becomes almost discontinuous. Thus, the system 
      of basis functions 
      should be enriched to account for this fact. The appropriate 
      choice of functions for such enrichment is a topic for further 
      research.
\item A rigorous proof of error estimates and the convergence to the
      Boltz\-mann-Maxwell equilibrium for the case of hard potentials with
      cut-off collision kernels is also a subject of an upcoming study.
\end{enumerate}
\section{Acknowledgement}
This research was supported by the the Institute for 
Computational Engineering and Sciences (ICES) the University of Texas at Austin  
through two research visits of the second author to ICES.
\noindent
We thank Dr. R. Grzhibovskis, Saarland University, for comments that greatly 
improved the manuscript.
\newpage
\bibliographystyle{plain}
\bibliography{Boltzmann}
\end{document}